# Expected extremal area of facets of random polytopes

Brett Leroux[*] Luis Rademacher[†] Carsten Schütt[‡]
Elisabeth M. Werner[§¶‖**]

**Abstract**

We study extremal properties of spherical random polytopes, the convex hull of random points chosen from the unit Euclidean sphere in $\mathbb{R}^n$. The extremal properties of interest are the expected values of the maximum and minimum surface area among facets. We determine the asymptotic growth in every fixed dimension, up to absolute constants.

# 1 Introduction

The research on random polytopes in convex bodies has been initiated in the works of Rényi and Sulanke [24, 25]. They studied the expected volume of the convex hull of finitely many points chosen from a convex body in dimension 2.

In higher dimensions this research has been extended by many people and by now the literature on random polytopes in convex bodies is vast. Schneider and Wieacker [26], Müller [19], Reitzner [22], Bárány [1, 2], and Schütt [27] extended the results of Rényi and Sulanke to convex bodies with smooth boundaries. In particular, it has been shown in [1, 2] and [27] that the expectation of the symmetric difference of a smooth convex body and a random polytope of $N$ points chosen from a convex body is of the order $N^{-\frac{2}{n-1}}$.

Bárány and Buchta settled the case of random polytopes in a polytope [3]. Bárány and Larman found a connection between random polytopes and floating bodies [4].

Much research focused on random polytopes whose vertices have been chosen from the interior of a convex body. Reitzner [21] and Schütt and Werner [28] studied the case of random polytopes whose vertices have been chosen from the boundary of a smooth convex

[*]Department of Mathematics, University of California, Davis. lerouxbew@gmail.com

[†]Department of Mathematics, University of California, Davis. lrademac@ucdavis.edu

[‡]Department of Mathematics, University of Kiel, Heinrich-Hecht-Platz 6, 24118 Kiel, Germany, schuett@math.uni-kiel.de

[§]Supported by NSF grant DMS-2103482

[¶]Department of Mathematics, Case Western Reserve University, 2145 Adalbert Road, Cleveland, OH 44106, USA, elisabeth.werner@case.edu

[‖]Work on this paper has been done while all authors visited ICERM, Providence, in the Fall of 2022. Also, work on this paper has been done while the last two authors visited the Hausdorff Institute, Bonn, in the Spring 2024.

[**]The first two authors were supported by the National Science Foundation under Grant CCF-2006994.

body. Then, amazingly, random approximation is almost as good as best approximation. In [17], Ludwig, Schütt and Werner approximated the Euclidean ball by arbitrarily positioned (neither inscribed nor circumscribed) random polytopes. This was generalized to smooth convex bodies by Grote and Werner in [15]. Random polytopes of vertices that are chosen from the boundary of a simple polytope have been investigated by Reitzner, Schütt and Werner [23]. Very recently Besau, Gusakova and Thäle [6] treated the case when the approximated convex body is neither smooth nor a polytope. Böröczky and Reitzner considered random polytopes that circumscribe a convex body [9].

Most results involve the symmetric difference metric, but there are also results where the error of approximation is measured in the Hausdorff distance. For instance Glasauer and Schneider [14] gave precise asymptotics for a smooth, convex body and random polytopes with respect to the Hausdorff distance. For polytopes in $\mathbb{R}^2$ such a study has been carried out by Prochno, Sonnleitner, Schütt and Werner [20]. Moreover there are results by Böröczky, Fodor, and Hug [10] when the error is measured in intrinsic volumes. Extensions to approximation of convex bodies on the sphere by random polytopes were obtained e.g., by Besau, Gusakova, Reitzner, Schütt, Thäle and Werner [5].

In this paper we consider extremal geometric and combinatorial properties of a random polytope on the sphere. A study of extremal properties of random polytopes has only been started recently [12]. Here we investigate the expected values of the maximum and minimum surface area among facets of a random polytope whose vertices are chosen randomly from the Euclidean sphere in $\mathbb{R}^n$. More precisely, we show that for a random polytope of $N$ chosen points the maximal expected surface area of a facet is of the order $\frac{\log N}{N}$ in dimension 2 and higher. Surprisingly, there is a change in behavior of the minimal expected surface area of a facet. It is of order $N^{-2}$ in dimension 2, of order $N^{-8/5}$ in dimension 3 and of order $N^{-3/2}$ in dimensions greater than or equal to 4.

**Our contributions.** We determine the asymptotic growth of the expected values of the maximum and minimum surface area among facets of a random polytope $P_N$, the convex hull of $N$ randomly chosen points on $S^{n-1}$, up to constants, in all dimensions. These results are stated in Theorems 1 and 2 below. There $\mathcal{F}(P)$ denotes the set of facets of polytope $P$.

**Theorem 1** (maximum area facet)**.** *Let $n \geq 2$. There exists constants $0 < a_n \leq b_n$ that depend only on $n$ such that for all $N \geq n+1$ we have*

$$a_n \frac{\log N}{N} \leq \mathbb{E} \max_{F \in \mathcal{F}(P_N)} \mathrm{vol}_{n-1}(F) \leq b_n \frac{\log N}{N}.$$

The proof of Theorem 1 is split into the upper bound and the lower bound, stated formally as Theorems 10 and 11. Theorem 10 (upper bound) is actually stronger in that it establishes a tail inequality on $\max_{F \in \mathcal{F}(P_N)} \mathrm{vol}_{n-1}(F)$.

**Theorem 2** (minimum area facet)**.** *There exist constants $0 < c_n \leq C_n$ and $N_n$ that depend only on $n$ such that:*

*(i) For $n = 2$ and all $N \geq 3$,*

$$3\sqrt{3}N^{-2} \leq \mathbb{E} \min_{F \in \mathcal{F}(P_N)} \mathrm{vol}_1(F) \leq 2\pi N^{-2}.$$

*(ii) For $n=3$ and all $N \geq N_3$,*

$$c_3 N^{-8/5} \leq \mathbb{E} \min_{F \in \mathcal{F}(P_N)} \mathrm{vol}_2(F) \leq C_3 N^{-8/5}.$$

*(iii) For $n \geq 4$ and all $N \geq N_n$,*

$$c_n N^{-3/2} \leq \mathbb{E} \min_{F \in \mathcal{F}(P_N)} \mathrm{vol}_{n-1}(F) \leq C_n N^{-3/2}.$$

We discuss those theorems. While the order of the expected minimal surface area of a facet is surprising, the order of the expected maximal facet was suggested by related, existing results. Indeed, Glasauer and Schneider [14] showed that the Hausdorff distance between the $n$-dimensional Euclidean unit ball $B_2^n$ and a random polytope $P_N$ of $N$ chosen points satisfies

$$\lim_{N\to\infty} \frac{d_H(B_2^n, P_N)}{(\frac{\log N}{N})^{\frac{2}{n-1}}} = \frac{1}{2}\left(\frac{2\sqrt{\pi}\Gamma(\frac{n+1}{2})}{\Gamma(\frac{n}{2})}\right)^{\frac{2}{n-1}},$$

where the limit is taken in the distributional sense. Therefore, for large $N$ there is a facet whose distance to the boundary is of the order $(\frac{\log N}{N})^{\frac{2}{n-1}}$. Such a facet is contained in a $(n-1)$-dimensional Euclidean ball with radius $\sqrt{2}(\frac{\log N}{N})^{\frac{1}{n-1}}$. Therefore the surface area of the facet is less than $2^{\frac{n-1}{2}} \mathrm{vol}_{n-1}(B_2^{n-1})\frac{\log N}{N}$. This suggests that the expected maximal facet has a surface area less than $2^{\frac{n-1}{2}} \mathrm{vol}_{n-1}(B_2^{n-1})\frac{\log N}{N}$.

One arrives at a similar conclusion by a result of Bonnet and O'Reilly [8] on convex hulls of random points from the sphere. They studied the expected number of facets whose signed distances of their affine hull to the origin are between two prescribed values $h_1$ and $h_2$. In particular, they showed that for properly chosen constants $c_1$ and $c_2$ these distances are almost all between

$$h_1 = \sqrt{1-\left(\frac{c_1 n}{N}\left(\log(N/n)\right)^{3/2}\right)^{\frac{2}{n-1}}} \qquad \text{and} \qquad h_2 = \sqrt{1-\left(\frac{c_2 n}{N}\right)^{2(n+1)/(n-1)^2}}.$$

Informally, we refer to $h_1$ as the minimal signed distance and to $h_2$ as the maximal signed distance.

This suggests that a facet with maximal surface area can be found at the minimal distance $h_1$ and a facet with minimal surface area can be found at the maximal distance $h_2$.

A facet at minimal distance $h_1$ is contained in an $(n-1)$-dimensional Euclidean ball of radius $\sqrt{1-h_1^2}$ and a facet at maximal distance $h_2$ is contained in an $(n-1)$-dimensional Euclidean ball with radius $\sqrt{1-h_2^2}$. This means that the facet at minimal distance has a surface area less than

$$\mathrm{vol}_{n-1}(B_2^{n-1})\frac{c_1 n(\log\frac{N}{n})^{3/2}}{N}$$

and a facet at maximal distance has at most area

$$\mathrm{vol}_{n-1}(B_2^{n-1})\left(\frac{c_2 n}{N}\right)^{\frac{n+1}{n-1}}.$$

The question arises whether the first expression gives the order of the maximal surface area and the second the minimal surface area. Neither is true. In the case of the maximal facet, the expression is close and only misses a factor of $\sqrt{\log N}$. In the case of the minimal facet the true order is very different. This means that the size of the minimal facet is not directly related to the maximal distance of the facet from the origin. Here we have a totally different phenomenon.

Note also that the behavior of the minimum is more intricate than that of the maximum: The asymptotic growth for the minimum depends on the dimension $n$.

**Outline of the paper.** Section 2 introduces notation and collects lemmas that will be used later. Sections 3 and 4 present the upper and lower bounds on the expected maximum facet volume, respectively. Section 5 gives an outline of the argument for the expected minimum facet volume. In Sections 6 to 9 we give the proofs of the upper and lower bounds on the expected minimum facet volume.

# 2 Preliminaries

Let $\chi(\cdot)$ be the indicator random variable of an event. We denote the Euclidean norm on $\mathbb{R}^n$ by $\|\cdot\|$. The convex hull of points $x_1, \ldots, x_N \in \mathbb{R}^n$ is $[x_1, \ldots, x_N]$. If $x_1, \ldots, x_n$ are linearly independent $\operatorname{aff}[x_1, \ldots, x_n]$ denotes the hyperplane containing the points $x_1, \ldots, x_n$. The Euclidean ball in $\mathbb{R}^n$ with center $x$ and radius $\rho$ is denoted by $B_2^n(x, \rho)$ and $B_2^n = B_2^n(0, 1)$. We write $\partial B_2^n$ or $S^{n-1}$ for the unit sphere. Let $\operatorname{vol}_n(K)$ be the $n$-dimensional volume of $K$. In long formulas we often write $|K|$ instead of $\operatorname{vol}_n(K)$. We have $\operatorname{vol}_n(B_2^n) = \pi^{n/2}/\Gamma(n/2)$ and $\operatorname{vol}_{n-1}(S^{n-1}) = n\pi^{n/2}/\Gamma(n/2)$. The hyperplane orthogonal to $\xi \neq 0$ with distance $p \geq 0$ from the origin is denoted by $H(\xi, p) = \{x : \langle \xi, x \rangle = p\|\xi\|\}$.

The Hausdorff distance between two convex bodies $C, K \subseteq \mathbb{R}^n$ is

$$d_H(C, K) = \inf\{\rho \geq 0 : C \subseteq K + \rho B_2^n \text{ and } K \subseteq C + \rho B_2^n\}.$$

When the ambient dimension $n$ is clear, $P_N$ denotes a random polytope defined as the convex hull of $N$ i.i.d. uniformly random points from $S^{n-1}$.

The set of all facets of a polytope $P$ is denoted by $\mathcal{F}(P)$.

**Gamma and Beta functions.** We will use the following properties of the Gamma and Beta functions: $\Gamma(x+1) = \int_0^\infty t^x e^{-t}\, \mathrm{d}t$, $\Gamma(n+1) = n!$ and $B(x, y) = \int_0^1 t^{x-1}(1-t)^{y-1}\, \mathrm{d}t = \frac{\Gamma(x)\Gamma(y)}{\Gamma(x+y)}$. We have $\lim_{x\to\infty} \frac{\Gamma(x+\alpha)}{x^\alpha \Gamma(x)} = 1$.

**Caps and halfspaces.** A *cap* of $S^{n-1}$ is the intersection of a closed half space with $S^{n-1}$. Equivalently, a cap is an intersection of $S^{n-1}$ with a Euclidean ball.

If $H$ is an affine hyperplane that does not contain the origin, then $H^-$ denotes the closed half space with $H$ as boundary that does not contain the origin. Similarly, $H^+$ denotes the closed half space with $H$ as boundary and containing the origin. If $H$ contains the origin then $H^+$ and $H^-$ refer to any of the two halfspaces with $H$ as boundary (the specific choice will not matter in our calculations because the event that $H$ contains 0 has measure 0 or because the calculation is indifferent to the choice).

By the *angle of a cap* we mean the polar angle, namely the angle between the rays from the center of the sphere to the apex of the cap (the pole) and the boundary sphere forming the base of the cap.

The *base* of a cap $S^{n-1} \cap H^-$ is $B_2^n \cap H$. The radius of the base is also referred to as the *radius* of the cap. The *height* of a cap is the perpendicular distance from the hyperplane that cuts the sphere to create the cap, to the topmost point of the cap. For a cap contained in a hemisphere, this means the height of a cap equals $1-p$, where $p$ is the distance of the hyperplane $H$ to the origin.

The next lemma is well known. We include a proof for completeness.

**Lemma 3.** *Let $n \geq 2$ and $0 \leq p \leq 1$. Let $H$ be a hyperplane in $\mathbb{R}^n$ with distance $p$ from the origin. The surface area of the cap $S^{n-1} \cap H^-$ is*

$$\mathrm{vol}_{n-1}(S^{n-1} \cap H^-) = \mathrm{vol}_{n-2}(S^{n-2}) \int_p^1 (1-s^2)^{\frac{n-3}{2}} \, \mathrm{d}s. \tag{1}$$

*Moreover, for $n \geq 4$,*

$$(1-p^2)^{\frac{n-1}{2}} |B_2^{n-1}| \leq \mathrm{vol}_{n-1}(S^{n-1} \cap H^-) \leq \frac{1}{p}(1-p^2)^{\frac{n-1}{2}} |B_2^{n-1}|. \tag{2}$$

*In dimension 3, $\mathrm{vol}_2(S^2 \cap H^-) = 2\pi(1-p)$ and in dimension 2, (1) equals $2 \arccos p$.*

*Proof.* Let $\theta \in S^{n-1}$ be the vector orthogonal to $H$. Let $\mathrm{Pr} : S^{n-1} \cap H^- \to H$ be the orthogonal projection. Since $0 \leq p$, the map Pr is a bijection. Then

$$\mathrm{vol}_{n-1}(S^{n-1} \cap H^-) = \int_{B_2^n \cap H} \frac{\mathrm{d}x}{\langle N(\mathrm{Pr}^{-1}(x)), \theta \rangle}, \tag{3}$$

where $N(\mathrm{Pr}^{-1}(x))$ is the normal at $\mathrm{Pr}^{-1}(x)$. For $x \in B_2^n \cap H$ let $r = \|x - \mathrm{Pr}(\theta)\|$ be the distance of $x$ to the center of the Euclidean ball $B_2^n \cap H$. We have $\langle N(\mathrm{Pr}^{-1}(x)), \theta \rangle = \sqrt{1-r^2}$ and, by passing to polar coordinates,

$$\mathrm{vol}_{n-1}(S^{n-1} \cap H^-) = \mathrm{vol}_{n-2}(S^{n-2}) \int_0^{\sqrt{1-p^2}} \frac{r^{n-2}}{\sqrt{1-r^2}} \, \mathrm{d}r.$$

Equation (1) follows with the substitution $s = \sqrt{1-r^2}$.

Equation (2) follows from (3). Indeed, $p \leq \langle N(\mathrm{Pr}^{-1}(x)), \theta \rangle \leq 1$ and consequently

$$\begin{aligned}(1-p^2)^{\frac{n-1}{2}} \mathrm{vol}_{n-1}(B_2^{n-1}) &= \int_{B_2^n \cap H} \mathrm{d}x \leq \int_{B_2^n \cap H} \frac{\mathrm{d}x}{\langle N(\mathrm{Pr}^{-1}(x)), \theta \rangle} \\ &\leq \frac{1}{p} \int_{B_2^n \cap H} \mathrm{d}x = \frac{1}{p}(1-p^2)^{\frac{n-1}{2}} \mathrm{vol}_{n-1}(B_2^{n-1}). \qquad \square\end{aligned}$$

**Lemma 4.** *Let $n \geq 2$ and $R \geq 2$. Let $C \subseteq S^{n-1}$ be a cap with angle $\phi$ and $|C| = |S^{n-1}|/R$. Then*

$$\left(\frac{1}{R} \cdot \frac{\mathrm{vol}_{n-1}(S^{n-1})}{\mathrm{vol}_{n-1}(B_2^{n-1})}\right)^{\frac{1}{n-1}} \leq \phi \leq \left(\frac{\pi}{2}\right)^{\frac{n-2}{n-1}} \left(\frac{1}{R} \cdot \frac{\mathrm{vol}_{n-1}(S^{n-1})}{\mathrm{vol}_{n-1}(B_2^{n-1})}\right)^{\frac{1}{n-1}}. \tag{4}$$

For $n = 2$ we have $\phi = \frac{\pi}{R}$.

*Proof.* Since the surface area of $C$ equals $\mathrm{vol}_{n-1}(S^{n-1})/R$ we get by Lemma 3

$$\frac{\mathrm{vol}_{n-1}(S^{n-1})}{R} = \mathrm{vol}_{n-2}(S^{n-2}) \int_{\cos\phi}^{1} (1 - s^2)^{\frac{n-3}{2}} \, \mathrm{d}s.$$

We put $s = \cos u$ and get $\frac{1}{R} \cdot \frac{\mathrm{vol}_{n-1}(S^{n-1})}{\mathrm{vol}_{n-2}(S^{n-2})} = \int_0^\phi \sin^{n-2} u \, \mathrm{d}u$. Since $2u/\pi \leq \sin u \leq u$ for $0 \leq u \leq \frac{\pi}{2}$

$$\frac{1}{n-1}\left(\frac{2}{\pi}\right)^{n-2} \phi^{n-1} = \left(\frac{2}{\pi}\right)^{n-2} \int_0^\phi u^{n-2} \, \mathrm{d}u \leq \int_0^\phi \sin^{n-2} u \, \mathrm{d}u \leq \int_0^\phi u^{n-2} \, \mathrm{d}u = \frac{1}{n-1}\phi^{n-1}$$

and we get

$$\left(\frac{n-1}{R} \cdot \frac{\mathrm{vol}_{n-1}(S^{n-1})}{\mathrm{vol}_{n-2}(S^{n-2})}\right)^{\frac{1}{n-1}} \leq \phi \leq \left(\frac{n-1}{R} \cdot \frac{\mathrm{vol}_{n-1}(S^{n-1})}{\mathrm{vol}_{n-2}(S^{n-2})}\right)^{\frac{1}{n-1}} \left(\frac{\pi}{2}\right)^{\frac{n-2}{n-1}}.$$

The claim follows. □

The following lemma on densest packing of the Euclidean sphere by caps is also well known. Again, we include its proof for completeness.

**Lemma 5.** *Let $R \geq 2$ and $n \geq 2$. Let $k \in \mathbb{N}$ be the maximal number such that there are caps $C_1, \dots, C_k$ of $S^{n-1}$ with*

$$\mathrm{int}(C_i) \cap \mathrm{int}(C_j) = \emptyset \qquad i \neq j \tag{5}$$

*and for all $i = 1, \dots, k$*

$$\mathrm{vol}_{n-1}(C_i) = |S^{n-1}|/R. \tag{6}$$

*Then, for sufficiently big $R$,*

$$3^{-n} R \leq k \leq R. \tag{7}$$

*Proof.* We show the right hand side inequality. We have $\bigcup_{i=1}^k C_i \subseteq S^{n-1}$. By (5) and (6), $k|S^{n-1}|/R \leq |S^{n-1}|$, which implies $k \leq R$.

Now we show the left hand side inequality. We construct caps $C_1, \dots, C_m$ satisfying (5), (6) and such that $3^{-n} R \leq m$. A maximal $\delta$-net in $S^{n-1}$ is a set of points $x_1, \dots, x_m \in S^{n-1}$ such that $\|x_i - x_j\| \geq \delta$ and for all $x \in S^{n-1}$ there is $i$, $1 \leq i \leq m$, such that $\|x - x_i\| \leq \delta$. Then the set of caps $C_i = B_2^n(x_i, \delta/2) \cap S^{n-1}$ for $i = 1, \dots, m$ satisfies (5) and (6). It is left to show $3^{-n} R \leq m$. Since $S^{n-1} \subseteq \bigcup_{i=1}^m B_2^n(x_i, \delta) \cap S^{n-1}$ we have $|S^{n-1}| \leq m|B_2^n(x_1, \delta) \cap S^{n-1}|$. The radius of the base of the cap $B_2^n(x_1, \delta) \cap S^{n-1}$ is less than $\delta$. By Lemma 3, $|B_2^n(x_1, \delta) \cap S^{n-1}| \leq \frac{\delta^{n-1}}{\sqrt{1-\delta^2}} |B_2^{n-1}|$, i.e. $\frac{\sqrt{1-\delta^2}}{\delta^{n-1}} \frac{|S^{n-1}|}{|B_2^{n-1}|} \leq m$. We choose $\delta$ such

that $|B_2^n(x_1,\delta/2)\cap S^{n-1}| = |S^{n-1}|/R$. The radius of the base of the cap $B_2^n\left(x_1,\frac{\delta}{2}\right)\cap S^{n-1}$ equals $\frac{\delta}{2}\sqrt{1-\frac{\delta^2}{16}}$. By Lemma 3

$$\left(\frac{\delta}{2}\right)^{n-1}\left(1-\frac{\delta^2}{16}\right)^{\frac{n-1}{2}}|B_2^{n-1}| \le \left|B_2^n\left(x_1,\frac{\delta}{2}\right)\cap S^{n-1}\right| = \frac{|S^{n-1}|}{R}.$$

For large $R$ we have $\delta\le\frac{1}{2}$ and $R\le\left(\frac{3}{\delta}\right)^{n-1}\frac{|S^{n-1}|}{|B_2^{n-1}|}\le 3^n m$. But we know $m\le k$. $\square$

**Lemma 6** (Miles [18])**.** *Let $\xi_1,\dots,\xi_{n+1}$ be randomly chosen points from $S^{n-1}$. Then*

$$\mathbb{E}\,\mathrm{vol}_n([\xi_1,\dots,\xi_{n+1}]) = \frac{1}{n!}\frac{\Gamma(\frac{n^2+1}{2})}{\Gamma(\frac{n^2}{2})}\left(\frac{\Gamma(\frac{n}{2})}{\Gamma(\frac{n+1}{2})}\right)^n\frac{\Gamma(\frac{n}{2})}{\Gamma(\frac{1}{2})}. \tag{8}$$

**Lemma 7** (Spherical Blaschke-Petkantschin, see e.g. [7])**.** *Let $f:(S^{n-1})^n\to[0,\infty)$ be a measurable function. Then*

$$\begin{aligned}&\int_{S^{n-1}}\cdots\int_{S^{n-1}} f(\zeta_1,\dots,\zeta_n)\,\mathrm{d}\zeta_1\cdots\mathrm{d}\zeta_n \qquad (9)\\ &= (n-1)!\int_0^1\int_{S^{n-1}}\int_{(H(\theta,p)\cap S^{n-1})^n} f(\xi_1,\dots,\xi_n)\frac{\mathrm{vol}_{n-1}([\xi_1,\dots,\xi_n])}{(1-p^2)^{n/2}}\,\mathrm{d}\xi_1\cdots\mathrm{d}\xi_n\,\mathrm{d}\theta\,\mathrm{d}p,\end{aligned}$$

*where $p$ is the distance of $H(\theta,p)$ from the origin, $\mathrm{d}\theta$ and $\mathrm{d}\zeta_1,\dots,\mathrm{d}\zeta_n$ are the $(n-1)$-dimensional Hausdorff measures and $\mathrm{d}\xi_1,\dots,\mathrm{d}\xi_n$ are the $(n-2)$-dimensional Hausdorff measures.*

# 3 Expected maximum upper bound

In this section we establish that $\mathbb{E}\max_{F\in\mathcal{F}(P_N)}\mathrm{vol}_{n-1}(F) = O(\frac{\log N}{N})$ (Theorem 10). The idea of the argument is to first prove a tail inequality on $d_H(P_N,B_2^n)$ (Lemma 9) and then deduce from it a tail inequality on $\max_{F\in\mathcal{F}(P_N)}\mathrm{vol}_{n-1}(F)$. This tail inequality implies an upper bound on the expectation from the identity $\mathbb{E}(X)=\int_0^\infty\mathbb{P}(X\ge t)\,\mathrm{d}t$.

The following Lemma gives an upper bound on the volume of a simplex.

**Lemma 8.** *Let $n\ge 2$ and let $\xi_1,\dots,\xi_n$ be points in $S^{n-1}$ such that $[\xi_1,\dots,\xi_n]$ is $(n-1)$-dimensional. Let $H=\mathrm{aff}\{\xi_1,\dots,\xi_n\}$ and suppose $0\notin H^-$. Let $\Delta$ be the height of the cap $S^{n-1}\cap H^-$. Then*

$$\mathrm{vol}_{n-1}([\xi_1,\dots,\xi_n]) \le \left(\frac{2n}{n-1}\right)^{\frac{n-1}{2}}\frac{\sqrt{n}}{(n-1)!}\Delta^{\frac{n-1}{2}}.$$

*Proof.* We may assume that $[\xi_1,\dots,\xi_n]$ is a regular simplex. Indeed, among all simplices contained in $B_2^n\cap H$ an inscribed regular simplex has the greatest volume. This follows from John's characterization of ellipsoids of minimal volume containing a convex body.

Let $\rho$ be the radius of $B_2^n \cap H$. Then $\Delta = 1 - \sqrt{1-\rho^2}$ and $\rho = \sqrt{2\Delta - \Delta^2}$. The volume of an $(n-1)$-dimensional regular simplex with sidelength $s$ is $\frac{\sqrt{n}}{(n-1)!}\left(\frac{s}{\sqrt{2}}\right)^{n-1}$. Therefore

$$\begin{aligned}\mathrm{vol}_{n-1}([\xi_1,\dots,\xi_n]) &\le \frac{\sqrt{n}}{(n-1)!}\left(\frac{n}{n-1}\right)^{\frac{n-1}{2}} \rho^{n-1} \\ &= \frac{\sqrt{n}}{(n-1)!}\left(\frac{n}{n-1}\right)^{\frac{n-1}{2}} (2\Delta-\Delta^2)^{\frac{n-1}{2}} \le \frac{\sqrt{n}}{(n-1)!}\left(\frac{n}{n-1}\right)^{\frac{n-1}{2}} (2\Delta)^{\frac{n-1}{2}}. \qquad \square\end{aligned}$$

We need the following estimate of the tail probability of the Hausdorff distance of a random polytope to the ball. Similar results can be found in [14]. For completeness we include a proof of the statement we need. Our proof is based on unpublished work [29].

**Lemma 9.** *Let* $0 < \Delta \le 1$. *Then*

$$\mathbb{P}\big(d_H(P_N, B_2^n) \ge \Delta\big) \le \left(1 - (\Delta/3)^{\frac{n-1}{2}} \frac{\mathrm{vol}_{n-1}(B_2^{n-1})}{\mathrm{vol}_{n-1}(S^{n-1})}\right)^N \frac{1}{(\Delta/3)^{\frac{n-1}{2}}} \frac{\mathrm{vol}_{n-1}(S^{n-1})}{\mathrm{vol}_{n-1}(B_2^{n-1})}.$$

*Proof.* Let $x_1,\dots,x_k$ be a maximal $\rho$-net in $S^{n-1}$. We choose $\rho = \sqrt{\Delta/2}$. Let $\xi_1,\dots,\xi_N$ be a choice from $S^{n-1}$ with $d_H\big([\xi_1,\dots,\xi_N], B_2^n\big) \ge \Delta$. Then there is a $\xi \in S^{n-1}$ such that all of the $\xi_1,\dots,\xi_N$ are elements of $H^+(\xi, 1-\Delta)$. Moreover,

$$H^-(\xi, 1-\Delta) \cap S^{n-1} = B_2^n\left(\xi, \sqrt{2\Delta}\right) \cap S^{n-1}. \tag{10}$$

Since $x_1,\dots,x_k$ is a maximal $\rho$-net, there is $i_0$ with $\|x_{i_0} - \xi\| \le \rho$ and, by triangle inequality, $B_2^n(x_{i_0},\rho) \subseteq B_2^n(\xi, 2\rho)$. Thus, $B_2^n\left(x_{i_0}, \sqrt{\Delta/2}\right) \cap S^{n-1} \subseteq B_2^n\left(\xi, \sqrt{2\Delta}\right) \cap S^{n-1}$. By assumption, the interior of (10) does not contain any $\xi_1,\dots,\xi_N$. It follows that the interior of $B_2^n\left(x_{i_0}, \sqrt{\Delta/2}\right) \cap S^{n-1}$ does not contain any $\xi_1,\dots,\xi_N$. Therefore

$$\begin{aligned}&\big\{(\xi_1,\dots,\xi_N) : d_H([\xi_1,\dots,\xi_N], B_2^n) \ge \Delta\big\} \\ &\subseteq \bigcup_{j=1}^k \left\{(\xi_1,\dots,\xi_N) : (\forall i = 1,\dots,N)\, \xi_i \notin B_2^n\left(x_j, \sqrt{\Delta/2}\right) \cap S^{n-1}\right\}.\end{aligned}$$

By Lemma 3,

$$\begin{aligned}\left|\left(B_2^n\left(x_j, \sqrt{\Delta/2}\right) \cap S^{n-1}\right)\right| &= \left|\left(H^-(x_j, \Delta/4) \cap S^{n-1}\right)\right| \ge \left(\frac{\Delta}{2} - \frac{\Delta^2}{16}\right)^{\frac{n-1}{2}} |B_2^{n-1}| \\ &\ge (\Delta/3)^{\frac{n-1}{2}} |B_2^{n-1}|. \end{aligned} \tag{11}$$

Consequently,

$$\mathbb{P}\big(d_H([\xi_1,\dots,\xi_N], B_2^n) \ge \Delta\big) \le k\left(1 - (\Delta/3)^{\frac{n-1}{2}} \frac{\mathrm{vol}_{n-1}(B_2^{n-1})}{\mathrm{vol}_{n-1}(S^{n-1})}\right)^N.$$

Since $x_1, \dots, x_k$ is a maximal $\rho$-net, the sets $\operatorname{int}\big(B_2^n(x_i, \frac{\rho}{2})\big) \cap S^{n-1}$, $i = 1, \dots, k$, are pairwise disjoint. By 11,

$$k\left(\frac{\Delta}{3}\right)^{\frac{n-1}{2}} |B_2^{n-1}| \le k\left(\frac{\Delta}{2} - \frac{\Delta^2}{16}\right)^{\frac{n-1}{2}} |B_2^{n-1}| \le \sum_{i=1}^{k} \left| B_2^n\left(x_i, \sqrt{\Delta/2}\right) \cap S^{n-1} \right| \le |S^{n-1}|,$$

and consequently $k \le (3/\Delta)^{\frac{n-1}{2}} |S^{n-1}|/|B_2^{n-1}|$. The claim follows. $\square$

**Theorem 10.** *Let $\xi_1, \dots, \xi_N$ be points chosen randomly from $S^{n-1}$ with respect to the uniform measure. Then*

$$\mathbb{P}\left(\max_{F\in\mathcal{F}([\xi_1,\dots,\xi_N])} \operatorname{vol}_{n-1}(F) \ge t\right) \le \min\left\{1, \frac{\sqrt{n}}{t\,(n-1)!}\left(\frac{6n}{n-1}\right)^{\frac{n-1}{2}} \frac{\operatorname{vol}_{n-1}(S^{n-1})}{\operatorname{vol}_{n-1}(B_2^{n-1})} \right. \tag{12}$$
$$\left. \times \left(1 - t\,\frac{(n-1)!}{\sqrt{n}}\left(\frac{n-1}{6n}\right)^{\frac{n-1}{2}} \frac{\operatorname{vol}_{n-1}(B_2^{n-1})}{\operatorname{vol}_{n-1}(S^{n-1})}\right)^N \right\}.$$

*Moreover*

$$\mathbb{E} \max_{F\in\mathcal{F}([\xi_1,\dots,\xi_N])} \operatorname{vol}_{n-1}(F) \le 6^{n/2} \frac{\sqrt{e\cdot n}}{(n-1)!} \frac{\operatorname{vol}_{n-1}(S^{n-1})}{\operatorname{vol}_{n-1}(B_2^{n-1})} \frac{\log N}{N}. \tag{13}$$

*Proof.* Let $F$ be a facet of $[\xi_1, \dots, \xi_N]$ and let $\Delta_F$ be the height of the cap given by $F$. By Lemma 8,

$$\mathbb{P}\left(\max_{F\in\mathcal{F}([\xi_1,\dots,\xi_N])} \operatorname{vol}_{n-1}(F) \ge t\right) \le \mathbb{P}\left(\left(\frac{2n}{n-1}\right)^{\frac{n-1}{2}} \frac{\sqrt{n}}{(n-1)!} \max_{F\in\mathcal{F}([\xi_1,\dots,\xi_N])} \Delta_F^{\frac{n-1}{2}} \ge t\right).$$

Since $d_H(B_2^n, [\xi_1, \dots, \xi_N]) \ge \max_{F\in\mathcal{F}([\xi_1,\dots,\xi_N])} \Delta_F$, it follows

$$\mathbb{P}\left(\max_{F\in\mathcal{F}([\xi_1,\dots,\xi_N])} \operatorname{vol}_{n-1}(F) \ge t\right) \le \mathbb{P}\left(d_H(B_2^n, [\xi_1, \dots, \xi_N]) \ge \frac{n-1}{2n}\left((n-1)!t/\sqrt{n}\right)^{\frac{2}{n-1}}\right).$$

By Lemma 9 the previous expression is smaller than

$$\frac{1}{t}\frac{\sqrt{n}}{(n-1)!}\left(\frac{6n}{n-1}\right)^{\frac{n-1}{2}} \frac{\operatorname{vol}_{n-1}(S^{n-1})}{\operatorname{vol}_{n-1}(B_2^{n-1})}\left(1 - t\frac{(n-1)!}{\sqrt{n}}\left(\frac{n-1}{6n}\right)^{\frac{n-1}{2}} \frac{\operatorname{vol}_{n-1}(B_2^{n-1})}{\operatorname{vol}_{n-1}(S^{n-1})}\right)^N.$$

Since the probability is less than 1 we have shown (12). Moreover,

$$\begin{aligned}
\mathbb{E} \max_{F\in\mathcal{F}([\xi_1,\dots,\xi_N])} \operatorname{vol}_{n-1}(F) &= \int_0^\infty \mathbb{P}\left(\max_{F\in\mathcal{F}([\xi_1,\dots,\xi_N])} \operatorname{vol}_{n-1}(F) \ge t\right) \mathrm{d}t \\
&\le \int_0^\infty \min\left\{1, \frac{1}{t}\frac{\sqrt{n}}{(n-1)!}\left(\frac{6n}{n-1}\right)^{\frac{n-1}{2}} \frac{\operatorname{vol}_{n-1}(S^{n-1})}{\operatorname{vol}_{n-1}(B_2^{n-1})}\right. \\
&\qquad \left.\left(1 - t\frac{(n-1)!}{\sqrt{n}}\left(\frac{n-1}{6n}\right)^{\frac{n-1}{2}} \frac{\operatorname{vol}_{n-1}(B_2^{n-1})}{\operatorname{vol}_{n-1}(S^{n-1})}\right)^N \right\} \mathrm{d}t.
\end{aligned}$$

We put

$$s = t\,\frac{(n-1)!}{\sqrt{n}}\left(\frac{n-1}{6n}\right)^{\frac{n-1}{2}}\frac{\mathrm{vol}_{n-1}(B_2^{n-1})}{\mathrm{vol}_{n-1}(S^{n-1})}$$

and get

$$\begin{aligned}\mathbb{E}\max_{F\in\mathcal{F}([\xi_1,\ldots,\xi_N])}\mathrm{vol}_{n-1}(F) &\le \frac{\sqrt{n}}{(n-1)!}\left(\frac{6n}{n-1}\right)^{\frac{n-1}{2}}\frac{\mathrm{vol}_{n-1}(S^{n-1})}{\mathrm{vol}_{n-1}(B_2^{n-1})}\int_0^\infty \min\left\{1,\frac{(1-s)^N}{s}\right\}\mathrm{d}s\\ &\le \frac{\sqrt{n}}{(n-1)!}\left(\frac{6n}{n-1}\right)^{\frac{n-1}{2}}\frac{\mathrm{vol}_{n-1}(S^{n-1})}{\mathrm{vol}_{n-1}(B_2^{n-1})}\int_0^\infty \min\left\{1,\frac{e^{-Ns}}{s}\right\}\mathrm{d}s.\end{aligned}$$

Since $(1+\frac{1}{k})^k$, $k\in\mathbb{N}$, is a monotonely increasing sequence converging to $e$,

$$\mathbb{E}\max_{F\in\mathcal{F}([\xi_1,\ldots,\xi_N])}\mathrm{vol}_{n-1}(F) \le 6^{\frac{n-1}{2}}\frac{\sqrt{e\cdot n}}{(n-1)!}\frac{\mathrm{vol}_{n-1}(S^{n-1})}{\mathrm{vol}_{n-1}(B_2^{n-1})}\int_0^\infty \min\left\{1,\frac{e^{-Ns}}{s}\right\}\mathrm{d}s. \tag{14}$$

Let $s_0$ be such that $\frac{e^{-Ns_0}}{s_0}=1$. The function $\frac{e^{-Ns}}{s}$ is decreasing on $(0,\infty)$. Therefore the right hand side of (14) equals

$$6^{\frac{n-1}{2}}\frac{\sqrt{e\cdot n}}{(n-1)!}\frac{\mathrm{vol}_{n-1}(S^{n-1})}{\mathrm{vol}_{n-1}(B_2^{n-1})}\left(s_0+\int_{s_0}^\infty \frac{1}{s}e^{-Ns}\,\mathrm{d}s\right). \tag{15}$$

Since $\frac{1}{s}e^{-Ns}$ is decreasing on $(0,\infty)$ we get

$$\frac{\log\frac{N}{\log N}}{N}\le s_0\le\frac{\log N}{N}.$$

Therefore (15) is less than

$$6^{\frac{n-1}{2}}\frac{\sqrt{e\cdot n}}{(n-1)!}\frac{\mathrm{vol}_{n-1}(S^{n-1})}{\mathrm{vol}_{n-1}(B_2^{n-1})}\left(\frac{\log N}{N}+\int_{\frac{1}{N}\log\frac{N}{\log N}}^\infty \frac{1}{s}e^{-Ns}\,\mathrm{d}s\right).$$

We have

$$\int_{\frac{1}{N}\log\frac{N}{\log N}}^\infty \frac{1}{s}e^{-Ns}\,\mathrm{d}s \le \frac{N}{\log\frac{N}{\log N}}\int_{\frac{1}{N}\log\frac{N}{\log N}}^\infty e^{-Ns}\,\mathrm{d}s = \frac{\log N}{N\log\frac{N}{\log N}}.$$

Therefore, for $N$ sufficiently big we have (13). □

# 4 Expected maximum lower bound

In this section we establish that $\mathbb{E}\max_{F\in\mathcal{F}(P_N)}\mathrm{vol}_{n-1}(F)=\Omega\big(\frac{\log N}{N}\big)$. This result is stated formally in Theorem 11 below and is one side of the inequality in (Theorem 1). The proof proceeds as follows. We first choose a set of $\Theta\big(\frac{N}{\log N}\big)$ disjoint caps on $S^{n-1}$. We then prove (in the proof of Theorem 11) a lower bound on the probability that at least one of these caps contains a facet with volume $\Omega\big(\frac{\log N}{N}\big)$. The main idea behind the proof of this lower bound is a specific application of the "second moment method", i.e., the fact that for a positive random variable $X$,

$$\mathbb{P}(X > 0) \geq \frac{(\mathbb{E}\, X)^2}{\mathbb{E}(X^2)},$$

which follows from the Cauchy-Schwarz inequality. In our application of the above inequality, $X$ is the number of caps that contain a facet with volume $\Omega\big(\frac{\log N}{N}\big)$. The proof in the case $n = 2$ is simpler and is dealt with separately.

**Theorem 11.** *Let $\xi_1, \ldots, \xi_N$ be randomly chosen points from $S^{n-1}$ with respect to the uniform measure. Then for $n = 2$ and sufficiently big $N$,*

$$\mathbb{E} \max_{F \in \mathcal{F}([\xi_1, \ldots, \xi_N])} \mathrm{vol}_{n-1}(F) \geq \frac{1}{2\pi} \frac{\log N}{N}.$$

*For $n \geq 3$ and for sufficiently big $N$,*

$$\mathbb{E} \max_{F \in \mathcal{F}([\xi_1, \ldots, \xi_N])} \mathrm{vol}_{n-1}(F) \geq \frac{\mathrm{vol}_{n-1}(S^{n-1})}{\mathrm{vol}_{n-1}(B_2^{n-1})} \frac{V_{1,n}^3}{V_{2,n}} \frac{2^{n^2-n-7}}{3^{n+3}\, \pi^{2n^2-3n-1}} \frac{\log N}{N}.$$

Before completing the proof of Theorem 11 we need two technical lemmas. The first is an easy bound on the distribution of the volume of a random simplex and the second establishes asymptotic upper and lower bounds for the probability that a given cap of $S^{n-1}$ contains a "large" facet. Let $\xi_1, \ldots, \xi_{n+1}$ be points that are chosen randomly from $S^{n-1}$. Let

$$V_{k,n} := \mathbb{E}\big((\mathrm{vol}_n[\xi_1 \ldots, \xi_{n+1}])^k\big). \tag{16}$$

**Lemma 12.** *Let $\xi_1, \ldots, \xi_{n+1}$ be points that are chosen randomly from $S^{n-1}$. Then*

$$\mathbb{P}\left\{\mathrm{vol}_n[\xi_1, \ldots, \xi_{n+1}] \geq \frac{V_{1,n}}{2}\right\} \geq \frac{V_{1,n}^2}{4V_{2,n}}.$$

*Proof.* By the Cauchy-Schwarz inequality, we have for any non-negative random variable $Z$ that $\mathbb{P}\{Z > \frac{1}{2}\,\mathbb{E}\, Z\} \geq \frac{(\mathbb{E}\, Z)^2}{4\,\mathbb{E}(Z^2)}$. The claim follows. □

Let $C$ be a cap of $S^{n-1}$ and let $R \geq 1$. We define the constant

$$\mathrm{vol}_{R,n} := \frac{\mathrm{vol}_{n-1}(S^{n-1})}{\mathrm{vol}_{n-1}(B_2^{n-1})} \frac{V_{1,n}}{2R}. \tag{17}$$

Let $H_{R,C} : (S^{n-1})^n \to \{0, 1\}$ be such that $H_{R,C}(\xi_1, \ldots, \xi_n) = 1$ if the following holds:

$$\mathrm{aff}(\xi_1, \ldots, \xi_n) \cap S^{n-1} \subseteq C \qquad \text{and} \qquad \mathrm{vol}_{n-1}[\xi_1, \ldots, \xi_n] \geq \mathrm{vol}_{R,n}\,. \tag{18}$$

Otherwise $H_{R,C}(\xi_1, \ldots, \xi_n) = 0$.

We do not specify $R$ for now. Later we shall choose $R := \frac{N}{\log N}$.

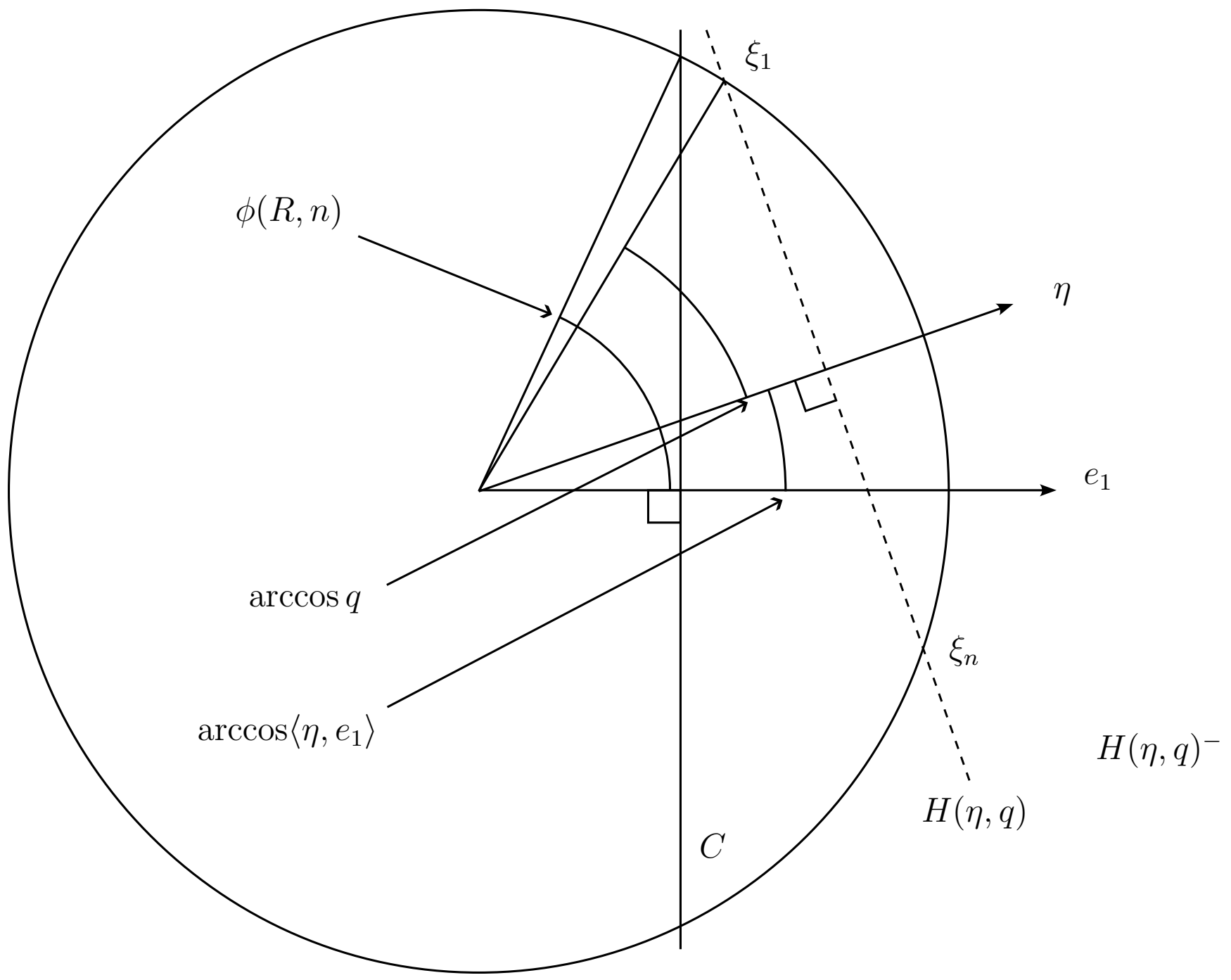

Figure 1: Angles in the proof of Lemma 13.

**Lemma 13.** *Let $n \geq 3$ and $R \geq \pi^n n$. Let $C$ be a cap of $S^{n-1}$ with $|C| = |S^{n-1}|/R$. Then*

$$\frac{1}{R^n}\frac{|S^{n-2}|V_{1,n}^3(n-1)!(n-1)^{n-3}}{3\cdot 2^5\cdot \pi^{n^2-2n-1}V_{2,n}} \leq \mathbb{E}\, H_{R,C} \leq \left(\frac{\pi}{2}\right)^{n(n-2)}(n-1)!(n-1)^{n-3}|S^{n-2}|V_{1,n}R^{-n}.$$

*Proof.* By definition

$$\mathbb{E}\, H_{R,C} = \frac{1}{|S^{n-1}|^n}\int_{S^{n-1}}\cdots\int_{S^{n-1}} H_{R,C}\,\mathrm{d}\xi_1\cdots\mathrm{d}\xi_n,$$

where $\mathrm{d}\xi_i$ denotes the restriction of the $(n-1)$-dimensional Hausdorff measure to $S^{n-1}$. We apply the spherical Blaschke-Petkantschin formula (9),

$$\mathbb{E}\, H_{R,C} = \frac{(n-1)!}{|S^{n-1}|^n}\int_0^1\int_{S^{n-1}}\int_{(H(\eta,q)\cap S^{n-1})^n} H_{R,C}(\xi_1,\dots,\xi_n)\frac{\big|[\xi_1,\dots,\xi_n]\big|}{(1-q^2)^{n/2}}\,\mathrm{d}\xi_1\cdots\mathrm{d}\xi_n\,\mathrm{d}\eta\,\mathrm{d}q, \quad (19)$$

where $\mathrm{d}\xi_i$ denotes the restriction of the $(n-2)$-dimensional Hausdorff measure to $H(\eta,q)\cap S^{n-1}$ and $\mathrm{d}\eta$ equals the restriction of the $(n-1)$-dimensional Hausdorff measure to $S^{n-1}$. We may assume that the hyperplane defining the cap $C$ is orthogonal to $e_1$. Let $\phi(R,n)$ be the angle of the cap $C$.

Let $\xi_1,\dots,\xi_n$ be affinely independent points from $S^{n-1}$ and let $H(\eta,q)$ be the hyperplane through $\xi_1,\dots,\xi_n$, i.e. the hyperplane orthogonal to unit vector $\eta$ at distance $q$ from the origin. Notice that $H_{R,C}(\xi_1,\dots,\xi_n)=1$ iff $\mathrm{vol}_{n-1}([\xi_1,\dots,\xi_n]) \geq \mathrm{vol}_{R,n}$ and the angle of the cap $S^{n-1} \cap H(\eta, q)^-$ (which equals $\arccos q$) is at most $\phi(R, n) - \arccos\langle \eta, e_1 \rangle$ (see Figure 1). The latter is equivalent to

$$\arccos\langle \eta, e_1 \rangle + \arccos q \leq \phi(R, n). \tag{20}$$

Therefore

$$H_{R,C}(\xi_1, \ldots, \xi_n) = \chi\big(\arccos\langle \eta, e_1 \rangle + \arccos q \leq \phi(R, n)\big)\, \chi\big(\mathrm{vol}_{n-1}([\xi_1, \ldots, \xi_n]) \geq \mathrm{vol}_{R,n}\big).$$

This in (19) gives

$$\begin{aligned} \mathbb{E}\, H_{R,C} = \frac{(n-1)!}{|S^{n-1}|^n} \int_0^1 \int_{S^{n-1}} \chi\big(\arccos\langle \eta, e_1 \rangle + \arccos q \leq \phi(N, n)\big) \qquad\qquad (21) \\ \times \int_{(H(\eta, q) \cap S^{n-1})^n} \chi\big(\mathrm{vol}_{n-1}([\xi_1, \ldots, \xi_n]) \geq \mathrm{vol}_{R,n}\big) \frac{\mathrm{vol}_{n-1}([\xi_1, \ldots, \xi_n])}{(1-q^2)^{n/2}}\, \mathrm{d}\xi_1 \cdots \mathrm{d}\xi_n\, \mathrm{d}\eta\, \mathrm{d}q. \end{aligned}$$

We first look at the **upper bound**: Since $\chi\big(\mathrm{vol}_{n-1}([\xi_1, \ldots, \xi_n]) \geq \mathrm{vol}_{R,n}\big) \leq 1$,

$$\begin{aligned} \mathbb{E}\, H_{R,C} \leq \frac{(n-1)!|S^{n-2}|^n}{|S^{n-1}|^n} \int_0^1 \int_{S^{n-1}} \chi\big(\arccos\langle \eta, e_1 \rangle + \arccos q \leq \phi(R, n)\big) \frac{(1-q^2)^{\frac{n^2-n-1}{2}}}{(1-q^2)^{n/2}} \\ \times J(q)\, \mathrm{d}\eta\, \mathrm{d}q, \end{aligned}$$

where

$$J(q) = \int_{(H(\eta, q) \cap S^{n-1})^n} \frac{\mathrm{vol}_{n-1}([\xi_1, \ldots, \xi_n])}{(1-q^2)^{\frac{n-1}{2}}} \frac{\mathrm{d}\xi_1}{|S^{n-2}|(1-q^2)^{\frac{n-2}{2}}} \cdots \frac{\mathrm{d}\xi_n}{|S^{n-2}|(1-q^2)^{\frac{n-2}{2}}}.$$

We have introduced the factor $(1-q^2)^{\frac{n^2-n-1}{2}}$ in order to normalize $J(q)$. Using the substitution $\xi_i = (1-q^2)^{\frac{1}{2}}\zeta_i$ and $\mathrm{d}\xi_i = (1-q^2)^{\frac{n-2}{2}}\, \mathrm{d}\zeta_i$, $i = 1, \ldots, n$, we get

$$J(q) = \frac{1}{|S^{n-2}|^n} \int_{S^{n-2}} \cdots \int_{S^{n-2}} \mathrm{vol}_{n-1}([\zeta_1, \ldots, \zeta_n])\, \mathrm{d}\zeta_1 \cdots \mathrm{d}\zeta_n = V_{1,n},$$

where $V_{1,n}$ is given by (16). Thus,

$$\begin{aligned} \mathbb{E}\, H_{R,C} \leq \frac{(n-1)!|S^{n-2}|^n}{|S^{n-1}|^n} V_{1,n} \\ \int_0^1 \int_{S^{n-1}} \chi\big(\arccos\langle \eta, e_1 \rangle + \arccos q \leq \phi(R, n)\big)(1-q^2)^{\frac{n^2-2n-1}{2}}\, \mathrm{d}\eta\, \mathrm{d}q. \end{aligned}$$

We use polar coordinates

$$\begin{aligned} \eta_1(\alpha_1, \ldots, \alpha_{n-1}) &= \cos\alpha_1 \\ \eta_k(\alpha_1, \ldots, \alpha_{n-1}) &= \sin\alpha_1 \cdots \sin\alpha_{k-1} \cos\alpha_k, \qquad k = 2, 3, \ldots, n-1 \\ \eta_n(\alpha_1, \ldots, \alpha_{n-1}) &= \sin\alpha_1 \cdots \sin\alpha_{n-1}. \end{aligned}$$

The change of coordinates induces the factor $\sin^{n-2}\alpha_1 \sin^{n-3}\alpha_2 \cdots \sin\alpha_{n-2}$. Then $\alpha_1 = \arccos\langle \eta, e_1\rangle$, i.e. $\cos\alpha_1 = \langle \eta, e_1\rangle$. We have

$$\begin{aligned}&\int_{S^{n-1}} \chi\big(\arccos\langle \eta, e_1\rangle + \arccos q \le \phi(R,n)\big)\,\mathrm{d}\eta \\ &= \int_0^{\pi} \cdots \int_0^{\pi}\int_0^{2\pi} \chi\big(\alpha_1 + \arccos q \le \phi(R,n)\big)\left(\prod_{k=1}^{n-2} \sin^{n-k-1}\alpha_k\right) \mathrm{d}\alpha_{n-1}\cdots \mathrm{d}\alpha_1.\end{aligned}$$

Since

$$|S^{n-2}| = \int_0^{\pi} \cdots \int_0^{\pi}\int_0^{2\pi} \left(\prod_{k=2}^{n-2} \sin^{n-k-1}\alpha_k\right) \mathrm{d}\alpha_{n-1}\cdots \mathrm{d}\alpha_2,$$

we get

$$\begin{aligned}&\int_{S^{n-1}} \chi\big(\arccos\langle \eta, e_1\rangle + \arccos q \le \phi(R,n)\big)\,\mathrm{d}\eta \\ &\qquad = |S^{n-2}| \int_0^{\pi} \chi\big(\alpha_1 + \arccos q \le \phi(R,n)\big) \sin^{n-2}\alpha_1\,\mathrm{d}\alpha_1.\end{aligned} \tag{22}$$

Therefore

$$\begin{aligned}\mathbb{E}\, H_{R,C} \le{}& \frac{(n-1)!|S^{n-2}|^{n+1}}{|S^{n-1}|^n} V_{1,n} \\ &\int_0^1 \int_0^{\pi} \chi\big(\alpha_1 + \arccos q \le \phi(R,n)\big)(1-q^2)^{\frac{n^2-2n-1}{2}} \sin^{n-2}\alpha_1\,\mathrm{d}\alpha_1\,\mathrm{d}q.\end{aligned}$$

Since $\chi\big(\alpha_1 + \arccos q \le \phi(R,n)\big) \le \chi\big(\alpha_1 \le \phi(R,n)\big)\,\chi\big(\arccos q \le \phi(R,n)\big)$, the double integral is less than or equal to

$$\int_0^{\pi} \chi\big(\alpha_1 \le \phi(R,n)\big) \sin^{n-2}\alpha_1\,\mathrm{d}\alpha_1 \times \int_0^1 \chi\big(\arccos q \le \phi(R,n)\big)(1-q^2)^{\frac{n^2-2n-1}{2}}\,\mathrm{d}q,$$

which, since $\phi(R,n) \le \pi$, is equal to

$$\int_0^{\phi(R,n)} \sin^{n-2}\alpha_1\,\mathrm{d}\alpha_1 \int_{\cos\phi(R,n)}^1 (1-q^2)^{\frac{n^2-2n-1}{2}}\,\mathrm{d}q.$$

We use twice that $\sin x \le x$ and substitute $q = \cos s$,

$$\begin{aligned}\mathbb{E}\, H_{R,C} &\le \frac{(n-1)!|S^{n-2}|^{n+1}}{|S^{n-1}|^n} V_{1,n} \left(\frac{\phi(R,n)^{n-1}}{n-1} \int_0^{\phi(R,n)} (\sin s)^{n^2-2n}\,\mathrm{d}s\right) \\ &\le \frac{(n-1)!|S^{n-2}|^{n+1}}{|S^{n-1}|^n} V_{1,n} \frac{\phi(R,n)^{n(n-1)}}{(n-1)^3}.\end{aligned}$$

Finally, this with (4) gives

$$\begin{aligned}\mathbb{E}\, H_{R,C} &\le \frac{(n-1)!|S^{n-2}|^{n+1}}{|S^{n-1}|^n} V_{1,n} \frac{1}{(n-1)^3} \left(\frac{\pi}{2}\right)^{n(n-2)} \left(\frac{1}{R}\,\frac{\mathrm{vol}_{n-1}(S^{n-1})}{\mathrm{vol}_{n-1}(B_2^{n-1})}\right)^n \\ &= (\pi/2)^{n(n-2)}\,(n-1)!(n-1)^{n-3}|S^{n-2}|V_{1,n}R^{-n}.\end{aligned}$$

Now we consider the **lower bound**: We start from (21). We restrict the integration with respect to $q$ to the interval $[0, f(N,n)]$,

$$\mathbb{E}\, H_{R,C} \geq \frac{(n-1)!|S^{n-2}|^n}{|S^{n-1}|^n} \int_0^{f(N,n)} \int_{S^{n-1}} \chi\big(\arccos\langle \eta, e_1\rangle + \arccos q \leq \phi(N,n)\big) \\ I(q)(1-q^2)^{\frac{n^2-2n-1}{2}}\, \mathrm{d}\eta\, \mathrm{d}q,$$

where

$$I(q) = \int_{(H(\eta,q)\cap S^{n-1})^n} \chi\big(\mathrm{vol}_{n-1}([\xi_1,\dots,\xi_n]) \geq \mathrm{vol}_{R,n}\big) \frac{\mathrm{vol}_{n-1}([\xi_1,\dots,\xi_n])}{(1-q^2)^{\frac{n-1}{2}}} \prod_{i=1}^n \frac{\mathrm{d}\xi_i}{|S^{n-2}|(1-q^2)^{\frac{n-2}{2}}}$$

and

$$f(N,n) := \sqrt{1-\left(\frac{1}{R}\frac{|S^{n-1}|}{|B_2^{n-1}|}\right)^{\frac{2}{n-1}}}. \tag{23}$$

We verify that $f(N,n)$ is well-defined, i.e. the expression under the square root is nonnegative. We have $|B_2^n| \leq \sqrt{\pi}|B_2^{n-1}|$. Indeed, since the Gamma function is increasing in $[3/2,\infty)$,

$$\frac{|B_2^n|}{|B_2^{n-1}|} = \frac{\pi^{\frac{n}{2}}\Gamma(\frac{n-1}{2}+1)}{\Gamma(\frac{n}{2}+1)\pi^{\frac{n-1}{2}}} \leq \sqrt{\pi}.$$

Moreover, since $R \geq n\cdot\pi^n$,

$$\left(\frac{1}{R}\frac{|S^{n-1}|}{|B_2^{n-1}|}\right)^{\frac{2}{n-1}} = \left(\frac{1}{R}\frac{n|B_2^n|}{|B_2^{n-1}|}\right)^{\frac{2}{n-1}} \leq \left(\frac{n\sqrt{\pi}}{R}\right)^{\frac{2}{n-1}} \leq \frac{1}{\pi^{2+\frac{1}{n-1}}}. \tag{24}$$

Using the substitution $\xi_i = (1-q^2)^{\frac{1}{2}}\zeta_i$ and $\mathrm{d}\xi_i = (1-q^2)^{\frac{n-2}{2}}\,\mathrm{d}\zeta_i$, $i=1,\dots,n$, we get

$$I(q) = \frac{1}{|S^{n-2}|^n} \int_{(S^{n-2})^n} \chi\left(\mathrm{vol}_{n-1}([\zeta_1,\dots,\zeta_n]) \geq \frac{\mathrm{vol}_{R,n}}{(1-q^2)^{\frac{n-1}{2}}}\right) \mathrm{vol}_{n-1}([\zeta_1,\dots,\zeta_n]) \prod_{i=1}^n \mathrm{d}\zeta_i.$$

Since $\chi\left(\mathrm{vol}_{n-1}([\zeta_1,\dots,\zeta_n]) \geq \mathrm{vol}_{R,n}(1-q^2)^{-\frac{n-1}{2}}\right)\mathrm{vol}_{n-1}([\zeta_1,\dots,\zeta_n])$ is nonnegative and the integral is an expectation,

$$\begin{aligned}
I(q) &= \int_0^\infty \mathbb{P}\left[\chi\left(\mathrm{vol}_{n-1}([\zeta_1,\dots,\zeta_n]) \geq \mathrm{vol}_{R,n}(1-q^2)^{-\frac{n-1}{2}}\right)\mathrm{vol}_{n-1}([\zeta_1,\dots,\zeta_n]) \geq s\right] \mathrm{d}s \\
&= \int_0^\infty \mathbb{P}\left[\mathrm{vol}_{n-1}([\zeta_1,\dots,\zeta_n]) \geq \max\left(s, \mathrm{vol}_{R,n}(1-q^2)^{-\frac{n-1}{2}}\right)\right] \mathrm{d}s \\
&\geq \int_0^{\mathrm{vol}_{R,n}(1-q^2)^{-\frac{n-1}{2}}} \mathbb{P}\left[\mathrm{vol}_{n-1}([\zeta_1,\dots,\zeta_n]) \geq \max\left(s, \mathrm{vol}_{R,n}(1-q^2)^{-\frac{n-1}{2}}\right)\right] \mathrm{d}s \\
&= \int_0^{\mathrm{vol}_{R,n}(1-q^2)^{-\frac{n-1}{2}}} \mathbb{P}\left[\mathrm{vol}_{n-1}([\zeta_1,\dots,\zeta_n]) \geq \mathrm{vol}_{R,n}(1-q^2)^{-\frac{n-1}{2}}\right] \mathrm{d}s.
\end{aligned}$$

If $q \le f(N,n)$, by (23) this means

$$(1-q^2)^{-\frac{n-1}{2}} \le \big(1-f(N,n)^2\big)^{-\frac{n-1}{2}} = R\,\mathrm{vol}_{n-1}(B_2^{n-1})/\mathrm{vol}_{n-1}(S^{n-1})$$

and so, using Lemma 12,

$$\begin{aligned}
&\mathbb{P}\big(\mathrm{vol}_{n-1}([\zeta_1,\dots,\zeta_n]) \ge \mathrm{vol}_{R,n}(1-q^2)^{-\frac{n-1}{2}}\big)\\
&\quad\ge \mathbb{P}\left(\mathrm{vol}_{n-1}([\zeta_1,\dots,\zeta_n]) \ge \mathrm{vol}_{R,n}\ R\,\mathrm{vol}_{n-1}(B_2^{n-1})/\mathrm{vol}_{n-1}(S^{n-1})\right)\\
&\quad= \mathbb{P}\big(\mathrm{vol}_{n-1}([\zeta_1,\dots,\zeta_n]) \ge (1/2)V_{1,n}\big) \ge (1/4)\,V_{1,n}^2/V_{2,n}.
\end{aligned}$$

It follows

$$I(q) \ge \frac{V_{1,n}^2\,\mathrm{vol}_{R,n}}{4V_{2,n}(1-q^2)^{\frac{n-1}{2}}} = \frac{V_{1,n}^3}{8V_{2,n}(1-q^2)^{\frac{n-1}{2}}}\frac{|S^{n-1}|}{|B_2^{n-1}|}\frac{1}{R}.$$

We have shown that

$$\begin{aligned}
\mathbb{E}\,H_{R,C} \ge{}& \frac{(n-1)(n-1)!\,|S^{n-2}|^{n-1}V_{1,n}^3}{8|S^{n-1}|^{n-1}V_{2,n}R}\\
&\int_0^{f(N,n)}\int_{S^{n-1}} \chi\big(\arccos\langle\eta,e_1\rangle + \arccos q \le \phi(R,n)\big)(1-q^2)^{\frac{n^2-3n}{2}}\,\mathrm{d}\eta\,\mathrm{d}q.
\end{aligned}$$

By (22) we get

$$\begin{aligned}
\mathbb{E}\,H_{R,C} \ge{}& \frac{(n-1)(n-1)!\,|S^{n-2}|^{n}V_{1,n}^3}{8|S^{n-1}|^{n-1}V_{2,n}R}\\
&\int_0^{f(N,n)}\int_0^{\pi} \chi\big(\alpha_1 + \arccos q \le \phi(R,n)\big)(\sin^{n-2}\alpha_1)(1-q^2)^{\frac{n^2-3n}{2}}\,\mathrm{d}\alpha_1\,\mathrm{d}q.
\end{aligned}$$

Since $\chi\big(\alpha_1 + \arccos q \le u\big) \ge \chi\big(\alpha_1 \le u/2\big)\,\chi\big(\arccos q \le u/2\big)$ we obtain

$$\begin{aligned}
\mathbb{E}\,H_{R,C} \ge{}& \frac{(n-1)(n-1)!\,|S^{n-2}|^{n}V_{1,n}^3}{8|S^{n-1}|^{n-1}V_{2,n}R}\int_0^{\pi}\chi\left(\alpha_1 \le \frac{\phi(R,n)}{2}\right)\sin^{n-2}\alpha_1\,\mathrm{d}\alpha_1\\
&\int_0^{f(N,n)}\chi\left(\arccos q \le \frac{\phi(R,n)}{2}\right)(1-q^2)^{\frac{n^2-3n}{2}}\,\mathrm{d}q.
\end{aligned}$$

Since $\phi(R,n) \le \pi$,

$$\mathbb{E}\,H_{R,C} \ge \frac{(n-1)(n-1)!\,|S^{n-2}|^{n}V_{1,n}^3}{8|S^{n-1}|^{n-1}V_{2,n}R}\int_0^{\phi(R,n)/2}\sin^{n-2}\alpha_1\,\mathrm{d}\alpha_1\int_{\cos\frac{\phi(R,n)}{2}}^{f(N,n)}(1-q^2)^{\frac{n^2-3n}{2}}\,\mathrm{d}q.$$

We use $\sin\alpha_1 \ge \frac{2}{\pi}\alpha_1$ for $0 \le \alpha_1 \le \frac{\pi}{2}$ and $q = \cos u$ to get

$$\begin{aligned}
\mathbb{E}\,H_{R,C} \ge{}& \frac{(n-1)(n-1)!\,|S^{n-2}|^{n}V_{1,n}^3}{8|S^{n-1}|^{n-1}V_{2,n}R}\left(\frac{2}{\pi}\right)^{n-2}\frac{1}{n-1}\left(\frac{\phi(R,n)}{2}\right)^{n-1}\\
&\int_{\arccos f(N,n)}^{\phi(R,n)/2}(\sin u)^{n^2-3n+1}\,\mathrm{d}u.
\end{aligned}$$

Thus

$$\mathbb{E}\, H_{R,C} \geq \frac{(n-1)!\,|S^{n-2}|^n V_{1,n}^3}{8|S^{n-1}|^{n-1} V_{2,n} R} \left(\frac{2}{\pi}\right)^{n^2-2n-1} \left(\frac{\phi(R,n)}{2}\right)^{n-1} \int_{\arccos f(N,n)}^{\phi(R,n)/2} u^{n^2-3n+1}\, \mathrm{d}u,$$

where

$$\int_{\arccos f(N,n)}^{\phi(R,n)/2} u^{n^2-3n+1}\, \mathrm{d}u = \frac{(\phi(R,n)/2)^{(n-1)(n-2)} - (\arccos f(N,n))^{(n-1)(n-2)}}{(n-1)(n-2)}.$$

We have $\arccos\sqrt{1-t} \leq \frac{\pi}{2}\sqrt{t}$ for $0 \leq t \leq 1$. Therefore

$$\arccos f(N,n) = \arccos\sqrt{1-(1-f(N,n)^2)} \leq \tfrac{\pi}{2}\sqrt{1-f(N,n)^2}.$$

By (23), $\arccos f(N,n) \leq \frac{\pi}{2} \cdot \left(\frac{1}{R}\frac{|S^{n-1}|}{|B_2^{n-1}|}\right)^{\frac{2}{n-1}}$. Therefore, using (4) to get a lower bound on $\phi(R,n)$,

$$\begin{aligned}
&\int_{\arccos f(N,n)}^{\phi(R,n)/2} u^{n^2-3n+1}\, \mathrm{d}u \\
&\geq \frac{1}{(n-1)(n-2)} \left(\frac{1}{2^{(n-1)(n-2)}} \left(\frac{1}{R} \cdot \frac{|S^{n-1}|}{|B_2^{n-1}|}\right)^{n-2} - \left(\frac{\pi}{2}\right)^{(n-1)(n-2)} \left(\frac{1}{R}\frac{|S^{n-1}|}{|B_2^{n-1}|}\right)^{2(n-2)}\right) \\
&= \frac{2^{-(n-1)(n-2)}}{(n-1)(n-2)} \left(\frac{1}{R} \cdot \frac{|S^{n-1}|}{|B_2^{n-1}|}\right)^{n-2} \left(1 - \left(\frac{\pi^{n-1}}{R}\frac{|S^{n-1}|}{|B_2^{n-1}|}\right)^{n-2}\right).
\end{aligned}$$

By (24) we get for $n \geq 3$

$$\pi^{(n-1)(n-2)} \cdot \left(\frac{1}{R}\frac{|S^{n-1}|}{|B_2^{n-1}|}\right)^{n-2} \leq \pi^{(n-1)(n-2)} \left(\frac{n\sqrt{\pi}}{R}\frac{\Gamma(\frac{n-1}{2}+1)}{\Gamma(\frac{n}{2}+1)}\right)^{n-2} \leq \pi^{-\frac{n-2}{2}}.$$

Therefore, since $\pi^{-\frac{1}{2}} \leq \frac{2}{3}$

$$\int_{\arccos f(N,n)}^{\phi(R,n)/2} u^{n^2-3n+1}\, \mathrm{d}u \geq \frac{2^{-(n-1)(n-2)}}{3(n-1)(n-2)} \left(\frac{1}{R} \cdot \frac{\mathrm{vol}_{n-1}(S^{n-1})}{\mathrm{vol}_{n-1}(B_2^{n-1})}\right)^{n-2}.$$

It follows

$$\begin{aligned}
\mathbb{E}\, H_{R,C} \geq{}& \frac{(n-1)(n-1)!\,|S^{n-2}|^n V_{1,n}^3}{8|S^{n-1}|^{n-1} V_{2,n} R} \left(\frac{2}{\pi}\right)^{n^2-2n-1} \frac{1}{n-1} 2^{-n+1} \frac{1}{R} \cdot \frac{\mathrm{vol}_{n-1}(S^{n-1})}{\mathrm{vol}_{n-1}(B_2^{n-1})} \\
&\times \frac{2^{-(n-1)(n-2)}}{3(n-1)(n-2)} \left(\frac{1}{R} \cdot \frac{\mathrm{vol}_{n-1}(S^{n-1})}{\mathrm{vol}_{n-1}(B_2^{n-1})}\right)^{n-2} \\
\geq{}& \frac{1}{R^n} \frac{|S^{n-2}| V_{1,n}^3 (n-1)!(n-1)^{n-3}}{3 \cdot 2^5 \cdot \pi^{n^2-2n-1} V_{2,n}}.
\end{aligned}$$

□

**Proof of Theorem 11**. For $n = 2$, this is essentially known and we deduce it from the expected length of the largest gap in a sample of uniformly random points in $S^1$. Let $X$ denote the (geodesic) length of the largest gap in a sample of $N$ uniformly random points from $S^1$. It is know that $\mathbb{E}(X) = \frac{1}{N}\sum_{i=1}^{N}\frac{1}{i}$. Let $Y$ denote the length of the longest edge of the convex hull of the same sample. Let $A$ be the event that the convex hull of the sample contains the origin. To relate the arc length with the edge length of the convex hull of the sample, notice that $\mathbb{P}(\bar{A}) = N/2^{N-1}$ (by Wendel's theorem). In general $Y$ can be much smaller than $X$ but if $A$ holds, then $Y \geq \frac{2}{\pi}X$. In this way, $\mathbb{E}(Y) \geq \mathbb{E}(Y \mid A)\mathbb{P}(A)$ with $\mathbb{P}(A) \geq 1/4$ for $N \geq 3$ while $\mathbb{E}(Y \mid A) \geq \frac{2}{\pi}\mathbb{E}(X \mid A)$. From the total expectation formula, $\mathbb{E}(X \mid A) = \frac{1}{\mathbb{P}(A)}(\mathbb{E}(X) - \mathbb{E}(X \mid \bar{A})\mathbb{P}(\bar{A})) \geq \mathbb{E}(X) - 2\pi\mathbb{P}(\bar{A})$. Combining everything, $\mathbb{E}(Y) \geq \frac{1}{2\pi}\mathbb{E}(X) - \frac{N}{2^{N-1}}$. Therefore $\mathbb{E}(Y) \geq \frac{1}{2\pi}\frac{(1/N)+\log(N)}{N} - \frac{N}{2^{N-1}}$. For sufficiently big $N$ we get $\mathbb{E}(Y) \geq \frac{1}{2\pi}\frac{\log(N)}{N}$.

The remainder of the proof is the case $n \geq 3$. Let $k \in \mathbb{N}$ be the maximal number such that there are caps $C_1, \dots, C_k$ of $S^{n-1}$ with $\mathrm{int}(C_i) \cap \mathrm{int}(C_j) = \emptyset,$, $i \neq j$, and for all $i = 1, \dots, k$ $|C_i| = |S^{n-1}|/R$. By Lemma 5,

$$3^{-n}R \leq k \leq R. \tag{25}$$

Let $G_i : S^{n-1} \times \cdots \times S^{n-1} \to \{0, 1\}$ such that $G_i(\xi_1, \dots, \xi_N) = 1$ if the following three conditions hold:

Exactly $n$ points $\xi_{\ell_1}, \dots, \xi_{\ell_n}$ are chosen from $C_i$ and they are affinely independent,

$$\mathrm{aff}(\xi_{\ell_1}, \dots, \xi_{\ell_n}) \cap S^{n-1} \subseteq C_i \qquad \text{and} \qquad \mathrm{vol}_{n-1}([\xi_{\ell_1}, \dots, \xi_{\ell_n}]) \geq \mathrm{vol}_{R,n}, \tag{26}$$

where $\mathrm{vol}_{R,n}$ is defined by (17). Else we put $G_i(\xi_1, \dots, \xi_N) = 0$. Note that if $G_i = 1$ then there is a facet of $[\xi_1, \dots, \xi_N]$ in $C_i$ with volume at least $\mathrm{vol}_{R,n}$.

We have for all $i = 1, \dots, k$,

$$\mathbb{E}(G_i) = \mathbb{P}(G_i = 1) = \binom{N}{n}\left(1 - \frac{1}{R}\right)^{N-n}\mathbb{E}\, H_{R,C_i} \tag{27}$$

where $H_{R,C_i}$ is given by (18). Indeed, the first equality holds since $G_i$ takes only the values 0 and 1. We have $G_i(\xi_1, \dots, \xi_N) = 1$ if and only if there are exactly $\xi_{\ell_1}, \dots, \xi_{\ell_n} \in C_i$ such that for all $\ell \notin \{\ell_1, \dots, \ell_n\}$ we have $x_\ell \notin C_i$ and such that $H_{R,C_i}(\xi_{\ell_1}, \dots, \xi_{\ell_n}) = 1$. We have $\binom{N}{n}$ choices of $\ell_1, \dots, \ell_n \in \{1, \dots, N\}$. The requirement that for all $\ell \notin \{\ell_1, \dots, \ell_n\}$ we have $x_\ell \notin C_i$ gives rise to the factor $(1 - \frac{1}{R})^{N-n}$. The conditions (26) mean that $H_{R,C_i}(\xi_{\ell_1}, \dots, \xi_{\ell_n}) = 1$. Thus (27) holds.

For all $i, j = 1, \dots, k$ with $i \neq j$

$$\mathbb{E}(G_i \cdot G_j) = \mathbb{P}(G_i \cdot G_j = 1) = \binom{N}{n}\binom{N-n}{n}\left(1 - \frac{2}{R}\right)^{N-2n}(\mathbb{E}\, H_{R,C_1})^2. \tag{28}$$

Indeed, $G_i \cdot G_j$ takes only the values 0 and 1. Therefore, $\mathbb{E}(G_i \cdot G_j) = \mathbb{P}(G_i \cdot G_j = 1)$. We have $G_i \cdot G_j = 1$ if and only if $G_i = 1$ and $G_j = 1$. This holds if and only if there are two disjoint sets $\xi_{\ell_1}, \dots, \xi_{\ell_n} \in C_i$ and $\xi_{\ell_{n+1}}, \dots, \xi_{\ell_{2n}} \in C_j$ such that for all $\ell \notin \{\ell_1, \dots, \ell_{2n}\}$ we have $\xi_\ell \notin C_i \cup C_j$ and (26) holds for $\xi_{\ell_1}, \dots, \xi_{\ell_n}$ and $\xi_{\ell_{n+1}}, \dots, \xi_{\ell_{2n}}$. There are $\binom{N}{n}\binom{N-n}{n}$ choices of $\ell_1, \dots, \ell_n$ and $\ell_{n+1}, \dots, \ell_{2n}$. The condition $\xi_\ell \notin C_i \cup C_j$ for $\ell \notin \{\ell_1, \dots, \ell_{2n}\}$ gives rise to the factor $\left(1 - \frac{2}{R}\right)^{N-2n}$. Since (26) holds for $\xi_{\ell_1}, \dots, \xi_{\ell_n}$ and $\xi_{\ell_{n+1}}, \dots, \xi_{\ell_{2n}}$ we have $H_{R,C_i}(\xi_{\ell_1}, \dots, \xi_{\ell_n}) = H_{R,C_j}(\xi_{\ell_{n+1}}, \dots, \xi_{\ell_{2n}}) = 1$. Moreover, $\mathbb{E}\, H_{R,C_i} = \mathbb{E}\, H_{R,C_j}$ for all $i, j = 1, \dots, k$. Thus (28) holds.

By (27)

$$\mathbb{E}\left(\sum_{i=1}^{k} G_i\right) = \sum_{i=1}^{k} \mathbb{E} G_i = k\binom{N}{n}\left(1 - \frac{1}{R}\right)^{N-n} \mathbb{E}\, H_{R,C_1}. \tag{29}$$

By (27), (28) and $G_i^2 = G_i$,

$$\mathbb{E}\left(\left(\sum_{i=1}^{k} G_i\right)^2\right) = \sum_{i,j=1}^{k} \mathbb{E}(G_i G_j) = \sum_{i \neq j} \mathbb{E}(G_i G_j) + \sum_{i=1}^{k} \mathbb{E}(G_i^2) \tag{30}$$

$$= k(k-1)\binom{N}{n}\binom{N-n}{n}\left(1 - \frac{2}{R}\right)^{N-2n} (\mathbb{E}\, H_{R,C_1})^2 + k\binom{N}{n}\left(1 - \frac{1}{R}\right)^{N-n} \mathbb{E}\, H_{R,C_1}.$$

By the Cauchy-Schwarz inequality and as $\mathbb{E}\, H_{R,C_1} > 0$,

$$\mathbb{P}\left(\sum_{i=1}^{k} G_i > 0\right) \geq \frac{\left(\sum_{i=1}^{k} \mathbb{E}\, G_i\right)^2}{\mathbb{E}\left(\left(\sum_{i=1}^{k} G_i\right)^2\right)}. \tag{31}$$

By (29) and (30),

$$\mathbb{P}\left(\sum_{i=1}^{k} G_i > 0\right) \geq \frac{\left[k\binom{N}{n}\left(1 - \frac{1}{R}\right)^{N-n} \mathbb{E}\, H_{R,C_1}\right]^2}{k(k-1)\binom{N}{n}\binom{N-n}{n}\left(1 - \frac{2}{R}\right)^{N-2n} (\mathbb{E}\, H_{R,C_1})^2 + k\binom{N}{n}\left(1 - \frac{1}{R}\right)^{N-n} \mathbb{E}\, H_{R,C_1}}$$

$$= \frac{k\binom{N}{n}\left(1 - \frac{1}{R}\right)^{2N-2n} \mathbb{E}\, H_{R,C_1}}{(k-1)\binom{N-n}{n}\left(1 - \frac{2}{R}\right)^{N-2n} \mathbb{E}\, H_{R,C_1} + \left(1 - \frac{1}{R}\right)^{N-n}}. \tag{32}$$

We choose $R = \frac{N}{\log N}$. The sequence $N^2(1 - \frac{\log N}{N})^{2N-2n}$, $N \geq 3$, is decreasing and

$$\lim_{N \to \infty} N^2 \left(1 - \frac{1}{R}\right)^{2N-2n} = \lim_{N \to \infty} N^2 \left(1 - \frac{\log N}{N}\right)^{2N-2n} = 1.$$

Therefore the numerator of (32) satisfies

$$k\binom{N}{n}\left(1 - \frac{1}{R}\right)^{2N-2n} \mathbb{E}\, H_{R,C_1} \geq \frac{1}{3^n}\frac{N}{\log N}\binom{N}{n}\frac{1}{N^2}\frac{(\log N)^n}{N^n}\frac{|S^{n-2}| V_{1,n}^3 (n-1)!(n-1)^{n-3}}{3 \cdot 2^5 \cdot \pi^{n^2-2n-1} V_{2,n}}$$

$$= \frac{1}{3^{n+1}}\frac{(\log N)^{n-1}}{N}\binom{N}{n}\frac{1}{N^n}\frac{|S^{n-2}| V_{1,n}^3 (n-1)!(n-1)^{n-3}}{2^5 \cdot \pi^{n^2-2n-1} V_{2,n}}.$$

For sufficiently large $N$

$$\binom{N}{n}\frac{1}{N^n}(n-1)!(n-1)^{n-3} \geq \frac{2}{9}(n-1)^{n-4}.$$

It follows

$$k\binom{N}{n}\left(1-\frac{1}{R}\right)^{2N-2n}\mathbb{E}\, H_{R,C_1} \geq \frac{1}{3^{n+3}}\frac{(\log N)^{n-1}}{N}\frac{|S^{n-2}|V_{1,n}^3(n-1)^{n-4}}{2^4\cdot\pi^{n^2-2n-1}V_{2,n}}. \tag{33}$$

and we have estimated the numerator of (32) from below.

Now we estimate the denominator of (32) from above. Clearly, for large enough $N$,

$$\left(1-\frac{1}{R}\right)^{N-n} \leq \exp\left(-\frac{(N-n)\log N}{N}\right) = \frac{1}{N}\exp\left(\frac{n\log N}{N}\right) \leq \frac{2}{N}. \tag{34}$$

We have for sufficiently large $N$

$$\left(1-\frac{2}{R}\right)^{N-2n} = \left(1-\frac{2\log N}{N}\right)^{N-2n} \leq \exp\left(-\frac{N-2n}{N}\log(N^2)\right) \leq \frac{2}{N^2}.$$

By Lemma 13 and (25),

$$\begin{aligned}(k-1)&\binom{N-n}{n}\left(1-\frac{2}{R}\right)^{N-2n}\mathbb{E}\, H_{R,C_1}\\ &\leq \frac{N}{\log N}\frac{N^n}{n!}\frac{2}{N^2}\frac{(\log N)^n}{N^n}\left(\frac{\pi}{2}\right)^{n(n-2)}(n-1)!(n-1)^{n-3}|S^{n-2}|V_{1,n}\\ &= 2\frac{(\log N)^{n-1}}{N}\left(\frac{\pi}{2}\right)^{n^2-n}\frac{(n-1)^{n-3}}{n}|S^{n-2}|V_{1,n}.\end{aligned}$$

So we have shown that for sufficiently big $N$,

$$\begin{aligned}&\mathbb{P}\left(\sum_{i=1}^k G_i > 0\right) \geq\\ &\frac{1}{3^{n+3}}\frac{(\log N)^{n-1}}{N}\frac{|S^{n-2}|V_{1,n}^3(n-1)^{n-4}}{2^4\cdot\pi^{n^2-2n-1}V_{2,n}}\left(2\frac{(\log N)^{n-1}}{N}\left(\frac{\pi}{2}\right)^{n^2-n}\frac{(n-1)^{n-3}}{n}|S^{n-2}|V_{1,n}+\frac{2}{N}\right)^{-1}\\ &\geq \frac{2^{n^2-n-6}}{3^{n+3}\,\pi^{2n^2-3n-1}}\frac{V_{1,n}^2}{V_{2,n}}.\end{aligned}$$

We have

$$\mathbb{P}\left(\max_{F\in\mathcal{F}(P_N)}\mathrm{vol}_{n-1}(F)\geq \mathrm{vol}_{R,n}\right) \geq \mathbb{P}\left(\sum_{i=1}^k G_i > 0\right). \tag{35}$$

Indeed, if $\sum_{i=1}^k G_i > 0$ then there is $i_0$ with $G_{i_0} > 0$ and consequently $G_{i_0} = 1$. This means that on the set $\sum_{i=1}^k G_i > 0$ we have $\max_{F\in\mathcal{F}(P_N)}\mathrm{vol}_{n-1}(F)\geq \mathrm{vol}_{R,n}$. By (35),

$$\mathbb{E}\max_{F\in\mathcal{F}(P_N)}\mathrm{vol}_{n-1}(F)\geq \mathrm{vol}_{R,n}\,\mathbb{P}\left(\max_{F\in\mathcal{F}(P_N)}\mathrm{vol}_{n-1}(F)\geq\mathrm{vol}_{R,n}\right)\geq \mathrm{vol}_{R,n}\,\mathbb{P}\left(\sum_{i=1}^k G_i>0\right),$$

where $\mathrm{vol}_{R,n}$ is defined in equation (17). Thus

$$\mathbb{E}\max_{F\in\mathcal{F}(P_N)}\mathrm{vol}_{n-1}(F)\geq \frac{\log N}{N}\frac{\mathrm{vol}_{n-1}(S^{n-1})}{\mathrm{vol}_{n-1}(B_2^{n-1})}\frac{V_{1,n}^3}{V_{2,n}}\ \frac{2^{n^2-n-7}}{3^{n+3}\,\pi^{2n^2-3n-1}}. \qquad \square$$

# 5 Strategy for the proof of the expected minimum

We prove the upper and lower bounds of Theorem 2 separately in Sections 7 and 8. The $n = 2$ case is known so here we focus on $n > 2$. We start with an outline of the arguments.

Assume that the $n$ points $\xi_1, \dots, \xi_n$ form a facet of the random polytope $P_N$ in $S^{n-1}$. The volume of the facet depends on both the height of the cap defined by the affine hull $H$ of the $n$ points as well as the shape of the simplex inside $H$. On average, we expect the volume of the facet to depend (up to constants depending only on the dimension) only on the height of the cap because the expected volume of a random simplex in $H \cap S^{n-1}$ is equal up to constants to the maximum volume of a simplex contained in $H \cap S^{n-1}$.

However, when one instead considers the *minimum* volume of a facet, we show that the volume does not depend only on the height of the cap. It also has a strong dependence on the shape of the simplex formed by the $n$ points inside the cap. In order to make this dependence precise, we first establish Lemma 14 which gives an approximation of the distribution of $\mathrm{vol}_{n-1}([\xi_1, \dots, \xi_n])$ where $\xi_1, \dots, \xi_n$ are points that are chosen randomly from $S^{n-1}$ with respect to the uniform measure. We show that when $n = 2$, the CDF $\mathbb{P}\big(\mathrm{vol}_2([\xi_1, \dots, \xi_n]) \leq t\big)$ is equal up to constants to $t^{2/3}$; when $n \geq 3$, the CDF is equal up to constants to $t$. This fact explains why the asymptotics in Theorem 2 are the same for all $n \geq 4$ but differ in the $n = 3$ case. We give an overview of the proof strategy for the upper and lower bounds in Theorem 2 below. For simplicity we describe the proof for the case $n \geq 4$. The proof in the $n = 3$ case is similar.

**Upper bound in Theorem 2.** We start by choosing $k \sim N$ pairwise disjoint caps on $S^{n-1}$ such that the volume of each cap is $|S^{n-1}|/N$. The proof shows that one of these caps is likely to contain a facet with the desired volume. This is done by first defining the function $\tilde{G}_{i,t}$, $1 \leq i \leq k$. For each $i$, $\tilde{G}_{i,t}$ is the indicator function for the event that (roughly speaking) the i-th cap contains a facet of $P_N$ and that the volume of the facet is less than $t/N$. The precise definition of $\tilde{G}_{i,t}$ is in Equations (72) to (74). The proof shows that

$$\mathbb{P}\left(\sum_{i=1}^{k} \tilde{G}_{i,t} > 0\right) \geq 1 - \frac{c}{Nt^2}$$

for some constant $c$ depending only on $n$. This inequality is established, as in Section 4, using the second moment method, i.e., the Cauchy Schwarz inequality.

Therefore, the main difficulty in the proof of this result is to obtain a suitable upper bound on $\mathbb{E}\big((\sum_{i=1}^{k} \tilde{G}_{i,t})^2\big)$. This is accomplished in Section 7 after first obtaining bounds on the closely related function $\tilde{H}_{t,C}$ which is introduced in (58) and analyzed in Lemma 15.

**Lower bound in Theorem 2.** The proof of the lower bound is somewhat simpler than the proof of the upper bound. We show in Proposition 16 that the probability that there exists a facet with volume less than $t$ is at most some constant times $t^2 N^3$. Then the proof of the lower bound on the minimum volume facet follows by setting $t = cN^{-3/2}$ and using Markov's inequality.

# 6 Technical Lemmas for the expected minimum

The first lemma establishes upper and lower bounds on the CDF of $\mathrm{vol}_{n-1}([\xi_1,\dots,\xi_n])$. Miles [18] computed the expected volume of a random simplex in a Euclidean ball. He did not use the CDF explicitly.

**Lemma 14.** *Let $\xi_1,\dots,\xi_{n+1}$ be points that are chosen randomly from $S^{n-1}$ with respect to the uniform measure.*
*(i) For $n=2$,*

$$(2t)^{2/3}/\pi^{8/3} \le \mathbb{P}\big(\mathrm{vol}_2([\xi_1,\xi_2,\xi_3]) \le t\big) \le 342\cdot t^{2/3}, \tag{36}$$

*where the left hand side inequality holds for all $t$ with $0\le t\le \pi$ and the right hand side inequality for all $t$ with $0\le t$.*
*(ii) For all $n\ge 3$,*

$$a_n t \le \mathbb{P}\big(\mathrm{vol}_n([\xi_1,\dots,\xi_{n+1}]) \le t\big) \le b_n t, \tag{37}$$

*where $0\le t\le |B_2^n|$ for the left hand side inequality and $0\le t$ for the right hand side inequality. The constants are*

$$a_n = \frac{1}{|B_2^n|}\left(\frac{1}{\sqrt{e}\sqrt{n}}\frac{|S^{n-2}|^{n+1}}{n^2|S^{n-1}|^n}\frac{\Gamma(\frac{(n-1)^2+1}{2})}{\Gamma(\frac{(n-1)^2}{2})}\left(\frac{\Gamma(\frac{n-1}{2})}{\Gamma(\frac{n}{2})}\right)^{n-1}\frac{\Gamma(\frac{n-1}{2})}{\Gamma(\frac{1}{2})}\right) \tag{38}$$

*and*

$$b_n = \frac{n!|S^{n-2}|^{n+1}}{|S^{n-1}|^n}B\left(\frac{1}{2},\frac{n^2-3n+2}{2}\right). \tag{39}$$

*Proof.* **(i)** We show first the inequality from below and we consider first the case $t\le 1$. The point $\xi_1$ is chosen arbitrarily. The points $\xi_2,\xi_3$ are chosen from the cap $S^1\cap H^-(\xi_1,1-2^{-1/3}t^{2/3})$. Here we used that $t\le 1$. The height of the cap $S^1\cap H^-(\xi_1,1-2^{-1/3}t^{2/3})$ is $2^{-1/3}t^{2/3}$. By Lemma 3, the surface area of the cap $S^1\cap H^-(\xi_1,1-2^{-1/3}t^{2/3})$ equals $2\arccos\left(1-t^{\frac{2}{3}}/2^{\frac{1}{3}}\right)$ which is greater than $2^{4/3}t^{1/3}$. Therefore, the measure of all choices $(\xi_1,\xi_2,\xi_3)$ as in this paragraph is greater than $(2^{4/3}t^{1/3})^2/(2\pi)^2 = 2^{2/3}t^{2/3}/\pi^2$. On the other hand,

$$\mathrm{vol}_2([\xi_1,\xi_2,\xi_3]) = (1/2)\,\mathrm{vol}_1([\xi_2,\xi_3])\,\mathrm{d}\big(\xi_1,\mathrm{aff}(\xi_2,\xi_3)\big) \le t.$$

Therefore, for $t\le 1$,

$$\mathbb{P}\big(\mathrm{vol}_2([\xi_1,\xi_2,\xi_3]) \le t\big) \ge (2t)^{2/3}/\pi^2.$$

We extend this to all $t$ with $0\le t\le\pi$. Let $1\le t\le\pi$. Then

$$\mathbb{P}\big(\mathrm{vol}_2([\xi_1,\xi_2,\xi_3]) \le t\big) \ge \mathbb{P}\big(\mathrm{vol}_2([\xi_1,\xi_2,\xi_3]) \le 1\big) \ge 2^{2/3}/\pi^2 \ge (2t)^{2/3}/\pi^{8/3}.$$

Now we show the inequality from above. We may assume that $\xi_1=e_1=(1,0)$. For any choice of $\xi_2$ we compute the measure of all $\xi_3$ such that $\mathrm{vol}_2([\xi_1,\xi_2,\xi_3])\le t$. We parametrize $\xi_2=(\cos\phi,\sin\phi)$ with $\phi\in[0,\pi]$. Because of symmetry it is enough to consider the range $[0,\pi]$ instead of $[0,2\pi]$. Then $\|\xi_1-\xi_2\| = 2\sin(\phi/2)$ and

$$t \ge \mathrm{vol}_2([\xi_1,\xi_2,\xi_3]) = (1/2)\|\xi_1-\xi_2\|\,\mathrm{d}\big(\xi_3,\mathrm{aff}(\xi_1,\xi_2)\big),$$

or

$$t/\sin(\phi/2) \geq \mathrm{d}\big(\xi_3, \mathrm{aff}(\xi_1,\xi_2)\big). \tag{40}$$

Since $\mathrm{d}\big(\xi_3, \mathrm{aff}(\xi_1,\xi_2)\big) \leq 2$, we have $\min\{2, t/\sin(\phi/2)\} \geq \mathrm{d}\big(\xi_3, \mathrm{aff}(\xi_1,\xi_2)\big)$. The hyperplane, i.e. the line through $\xi_1$ and $\xi_2$ is given by $H\left((\cos\frac{\phi}{2},\sin\frac{\phi}{2}),\cos\frac{\phi}{2}\right)$. By (40) we can choose any $\xi_3$ from the set

$$S^1 \cap H^-\left((\cos\tfrac{\phi}{2},\sin\tfrac{\phi}{2}),\cos\tfrac{\phi}{2}-\tfrac{t}{\sin\frac{\phi}{2}}\right)\cap H^+\left((\cos\tfrac{\phi}{2},\sin\tfrac{\phi}{2}),\cos\tfrac{\phi}{2}+\tfrac{t}{\sin\frac{\phi}{2}}\right), \tag{41}$$

where both half spaces contain the point $\cos\frac{\phi}{2}(\cos\frac{\phi}{2},\sin\frac{\phi}{2})$. It follows,

$$\begin{aligned}&\mathbb{P}\big(\mathrm{vol}_2([\xi_1,\xi_2,\xi_3])\leq t\big) \qquad (42)\\ &=\frac{1}{\pi}\int_0^{\pi}\left|S^1 \cap H^-\left((\cos\tfrac{\phi}{2},\sin\tfrac{\phi}{2}),\cos\tfrac{\phi}{2}-\tfrac{t}{\sin\frac{\phi}{2}}\right)\cap H^+\left((\cos\tfrac{\phi}{2},\sin\tfrac{\phi}{2}),\cos\tfrac{\phi}{2}+\tfrac{t}{\sin\frac{\phi}{2}}\right)\right|\mathrm{d}\phi.\end{aligned}$$

We split the integral into two parts, from 0 to $10\,\phi_0$ and $10\,\phi_0$ to $\pi$. We explain why we do this. If $\cos\frac{\phi}{2}+\frac{t}{\sin\frac{\phi}{2}}\geq 1$, then (41) equals

$$S^1 \cap H^-\left((\cos\tfrac{\phi}{2},\sin\tfrac{\phi}{2}),\cos\tfrac{\phi}{2}-\tfrac{t}{\sin(\phi/2)}\right). \tag{43}$$

Let $\phi_0$ be the unique solution of $t=\sin\frac{\phi}{2}(1-\cos\frac{\phi}{2})$. The value $\phi_0$ is the exact value such that the set (43) equals the set (41) for all $\phi$ with $0\leq\phi\leq\phi_0$. Since $\sin x\geq\frac{2x}{\pi}$ for $0\leq x\leq\frac{\pi}{2}$ and $\cos x\leq 1-\frac{x^2}{\pi}$ for $0\leq x\leq\frac{\pi}{2}$ this implies $t\geq\frac{\phi_0^3}{4\pi^2}$ (we apply the inequalities to $x=\phi/2$ so that we cover all $\phi\in[0,\pi]$). Since $\sin x\leq x$ for $x\geq 0$ and $\cos x\geq 1-\frac{x^2}{2}$ for $x\geq 0$ we get $t\leq\phi_0^3/16$. Altogether,

$$2^{4/3}t^{1/3}\leq\phi_0\leq(2\pi)^{2/3}t^{1/3}. \tag{44}$$

The cap (43) has height $1-\cos\frac{\phi}{2}+\frac{t}{\sin(\phi/2)}$. By Lemma 3, its surface area is equal to

$$2\arccos\left(\cos\frac{\phi}{2}-\frac{t}{\sin(\phi/2)}\right)\leq 2\arccos\left(1-\frac{\phi^2}{8}-\frac{\pi t}{\phi}\right),$$

using the fact that arccos is decreasing, $\cos x\geq 1-\frac{x^2}{2}$ and $\sin x\geq\frac{2x}{\pi}$. Since $\cos x\leq 1-\frac{x^2}{\pi}$ and arccos is decreasing we have $x\geq\arccos\big(1-\frac{x^2}{\pi}\big)$. Using this and that the square root is a concave function,

$$2\arccos\left(1-\frac{\phi^2}{8}-\frac{\pi t}{\phi}\right)\leq 2\sqrt{\frac{\pi}{8}\phi^2+\pi^2\frac{t}{\phi}}\leq\sqrt{\frac{\pi}{2}}\phi+2\pi\sqrt{\frac{t}{\phi}}.$$

By this and (44),

$$\begin{aligned}\int_0^{10\phi_0}2\arccos\left(1-\frac{\phi^2}{8}-\frac{\pi t}{\phi}\right)\mathrm{d}\phi &\leq \int_0^{10(2\pi)^{2/3}t^{1/3}}\sqrt{\frac{\pi}{2}}\phi+2\pi\sqrt{\frac{t}{\phi}}\,\mathrm{d}\phi \qquad (45)\\ &= t^{\frac{2}{3}}\pi^{\frac{11}{6}}\frac{100}{2^{\frac{1}{6}}}+t^{\frac{2}{3}}2^{\frac{17}{6}}\sqrt{5}\pi^{\frac{4}{3}}\leq 341\,\pi\cdot t^{\frac{2}{3}}.\end{aligned}$$

This leaves us with the integral

$$\int_{10\phi_0}^{\pi} 2\arccos\left(\cos\frac{\phi}{2} - \frac{t}{\sin(\phi/2)}\right) - 2\arccos\left(\cos\frac{\phi}{2} + \frac{t}{\sin(\phi/2)}\right) \mathrm{d}\phi. \tag{46}$$

By the mean value theorem there is $\alpha$ with $\cos\frac{\phi}{2} - \frac{t}{\sin(\phi/2)} \le \alpha \le \cos\frac{\phi}{2} + \frac{t}{\sin(\phi/2)}$ and

$$\frac{\arccos\left(\cos\frac{\phi}{2} + \frac{t}{\sin(\phi/2)}\right) - \arccos\left(\cos\frac{\phi}{2} - \frac{t}{\sin(\phi/2)}\right)}{2t/\sin(\phi/2)} = \arccos'(\alpha) = -\frac{1}{\sqrt{1-\alpha^2}}.$$

As $\frac{1}{\sqrt{1-\alpha^2}}$ is an increasing function on $[0,1)$, we get

$$\begin{aligned}\frac{\arccos\left(\cos\frac{\phi}{2} - \frac{t}{\sin(\phi/2)}\right) - \arccos\left(\cos\frac{\phi}{2} + \frac{t}{\sin(\phi/2)}\right)}{2t/\sin(\phi/2)} &= \frac{1}{\sqrt{1-\alpha^2}} \\ &\le \frac{1}{\sqrt{1-\left(\cos\frac{\phi}{2} + \frac{t}{\sin(\phi/2)}\right)^2}}.\end{aligned}$$

We have, using $\frac{2}{\pi}x \le \sin x$,

$$1 - \left(\cos\frac{\phi}{2} + \frac{t}{\sin(\phi/2)}\right)^2 = (\sin\tfrac{\phi}{2})^2 - 2t\frac{\cos\frac{\phi}{2}}{\sin\frac{\phi}{2}} - \frac{t^2}{(\sin\frac{\phi}{2})^2} \ge \frac{\phi^2}{\pi^2} - 2\pi\frac{t}{\phi} - \frac{t^2\pi^2}{\phi^2}.$$

By (44), $\phi \ge 10\phi_0 \ge 10\cdot 2^{4/3}t^{1/3}$, i.e. $t \le \phi^3/16000$. Therefore, the previous expression is greater than

$$\frac{\phi^2}{\pi^2} - \pi\frac{\phi^2}{8000} - \frac{\phi^4\pi^2}{16000^2}.$$

For $0 \le \phi \le \pi$ we have $\phi^4 \le \pi^2\phi^2$. Thus,

$$\frac{\phi^2}{\pi^2} - \pi\frac{\phi^2}{8000} - \frac{\phi^2\pi^4}{16000^2} = \phi^2\left(\frac{1}{\pi^2} - \frac{\pi}{8000} - \frac{\pi^4}{16000^2}\right) \ge \phi^2\left(\frac{1}{10} - \frac{1}{2000} - \frac{1}{1600^2}\right) \ge \frac{\phi^2}{20}$$

and therefore

$$\sqrt{1-\left(\cos\frac{\phi}{2} + \frac{t}{\sin(\phi/2)}\right)^2} \ge \frac{\phi}{5}.$$

Altogether,

$$\arccos\left(\cos\frac{\phi}{2} - \frac{t}{\sin(\phi/2)}\right) - \arccos\left(\cos\frac{\phi}{2} + \frac{t}{\sin(\phi/2)}\right) \le \frac{2t}{\frac{\phi}{5}\sin\frac{\phi}{2}} \le \frac{10\pi t}{\phi^2}.$$

As $t \le \pi$ and with (44), the integral (46) is at most

$$\int_{10\cdot 2^{4/3}t^{1/3}}^{\pi} \frac{20\pi t}{\phi^2}\,\mathrm{d}\phi = \frac{\pi}{2^{1/3}}t^{2/3} - 20t = \frac{\pi}{2^{1/3}}t^{2/3}\left(1 - \frac{20\cdot 2^{1/3}}{\pi}t^{1/3}\right).$$

For $t \leq 0.09$

$$2 \int_{10\phi_0}^{\pi} \arccos\left(\cos\tfrac{\phi}{2} - \frac{t}{\sin\frac{\phi}{2}}\right) - \arccos\left(\cos\tfrac{\phi}{2} + \frac{t}{\sin\frac{\phi}{2}}\right) \mathrm{d}\phi \leq \frac{\pi}{2^{\frac{4}{3}}} t^{\frac{2}{3}}. \tag{47}$$

By (42), (45) and (47)

$$\mathbb{P}\big(\mathrm{vol}_2([\xi_1, \xi_2, \xi_3]) \leq t\big) \leq 341 \cdot t^{\frac{2}{3}} + \frac{1}{2^{\frac{4}{3}}} t^{\frac{2}{3}} \leq 342 \cdot t^{\frac{2}{3}}.$$

Thus we have verified (36) for $0 \leq t \leq 0.09$. We extend this to all $t$ with $0 \leq t \leq \pi$. For $t$ with $t \geq 0.09$ we have

$$\mathbb{P}(\mathrm{vol}_2([\xi_1, \xi_2, \xi_3]) \leq t) \leq 1 \leq 38 \cdot t^{\frac{2}{3}}.$$

**(ii)** First we show the left hand side estimate of (37). Clearly,

$$\mathrm{vol}_n([\xi_1, \dots, \xi_{n+1}]) = (1/n)\, \mathrm{d}\big(\xi_{n+1}, \mathrm{aff}(\xi_1, \dots, \xi_n)\big)\, \mathrm{vol}_{n-1}([\xi_1, \dots, \xi_n])$$

and therefore

$$\mathbb{P}\big(\mathrm{vol}_n([\xi_1, \dots, \xi_{n+1}]) \leq t\big) = \mathbb{P}\left\{\mathrm{d}\big(\xi_{n+1}, \mathrm{aff}(\xi_1, \dots, \xi_n)\big) \leq \frac{nt}{\mathrm{vol}_{n-1}([\xi_1, \dots, \xi_n])}\right\}. \tag{48}$$

Trivially,

$$\mathbb{P}\big(\mathrm{vol}_n([\xi_1, \dots, \xi_{n+1}]) \leq t\big) \geq \mathbb{P}\left\{\mathrm{vol}_n([\xi_1, \dots, \xi_{n+1}]) \leq t \text{ and } \mathrm{d}\big(0, \mathrm{aff}(\xi_1, \dots, \xi_n)\big) \leq 1/n\right\}.$$

With (48)

$$\begin{aligned} &\mathbb{P}\big(\mathrm{vol}_n([\xi_1, \dots, \xi_{n+1}]) \leq t\big) \\ &\geq \mathbb{P}\left\{\mathrm{d}\big(\xi_{n+1}, \mathrm{aff}(\xi_1, \dots, \xi_n)\big) \leq \frac{nt}{\mathrm{vol}_{n-1}([\xi_1, \dots, \xi_n])} \,\middle|\, \mathrm{d}(0, \mathrm{aff}(\xi_1, \dots, \xi_n)) \leq \frac{1}{n}\right\} \\ &\quad \times \mathbb{P}\left\{\mathrm{d}\big(0, \mathrm{aff}(\xi_1, \dots, \xi_n)\big) \leq 1/n\right\}. \end{aligned} \tag{49}$$

By the Blaschke-Petkantschin formula (9),

$$\begin{aligned} &\mathbb{P}\left\{\mathrm{d}\big(0, \mathrm{aff}(\xi_1, \dots, \xi_n)\big) \leq 1/n\right\} \\ &= \frac{1}{|S^{n-1}|^n} \int_{S^{n-1}} \cdots \int_{S^{n-1}} \chi\left(\mathrm{d}\big(0, \mathrm{aff}(\xi_1, \dots, \xi_n)\big) \leq \frac{1}{n}\right) \mathrm{d}\xi_1 \cdots \mathrm{d}\xi_n \\ &= \frac{(n-1)!}{|S^{n-1}|^n} \int_0^{\frac{1}{n}} \int_{S^{n-1}} \int_{H(\theta,p)\cap S^{n-1}} \cdots \int_{H(\theta,p)\cap S^{n-1}} \frac{\mathrm{vol}_{n-1}([\xi_1, \dots, \xi_n])}{(1-p^2)^{\frac{n}{2}}} \mathrm{d}\xi_1 \cdots \mathrm{d}\xi_n \, \mathrm{d}\theta \, \mathrm{d}p. \end{aligned}$$

The $n$-fold integral over $H(\theta, p) \cap S^{n-1}$ is independent of $\theta$. Therefore the above expression equals

$$\frac{(n-1)!}{|S^{n-1}|^{n-1}} \int_0^{1/n} \int_{H(e_1,p)\cap S^{n-1}} \cdots \int_{H(e_1,p)\cap S^{n-1}} \frac{\mathrm{vol}_{n-1}([\xi_1, \dots, \xi_n])}{(1-p^2)^{n/2}} \mathrm{d}\xi_1 \cdots \mathrm{d}\xi_n \, \mathrm{d}p.$$

The set $H(e_1,p)\cap S^{n-1}$ is a Euclidean sphere with radius $\sqrt{1-p^2}$. We put $\xi_i = \eta_i\sqrt{1-p^2}$ and have $\mathrm{d}\xi_i = (1-p^2)^{\frac{n-2}{2}}\,\mathrm{d}\eta_i$. The above expression equals

$$\frac{(n-1)!}{|S^{n-1}|^{n-1}}\int_0^{1/n}\int_{S^{n-2}}\cdots\int_{S^{n-2}}\mathrm{vol}_{n-1}([\eta_1,\dots,\eta_n])(1-p^2)^{\frac{n^2-2n-1}{2}}\,\mathrm{d}\eta_1\cdots\mathrm{d}\eta_n\,\mathrm{d}p.$$

We apply Lemma 6 in dimension $n-1$ to get

$$\begin{aligned}&\mathbb{P}\left\{\mathrm{d}\big(0,\mathrm{aff}(\xi_1,\dots,\xi_n)\big)\le 1/n\right\}\\ &\quad=\frac{|S^{n-2}|^n}{|S^{n-1}|^{n-1}}\frac{\Gamma(\frac{n^2-2n+2}{2})}{\Gamma(\frac{(n-1)^2}{2})}\left(\frac{\Gamma(\frac{n-1}{2})}{\Gamma(\frac{n}{2})}\right)^{n-1}\frac{\Gamma(\frac{n-1}{2})}{\Gamma(\frac{1}{2})}\int_0^{\frac{1}{n}}(1-p^2)^{\frac{n^2-2n-1}{2}}\,\mathrm{d}p.\end{aligned}\tag{50}$$

By Bernoulli's inequality,

$$\frac{1}{2n}\le\frac{1}{n}\left(1-\frac{1}{n^2}\right)^{\frac{n^2-2n-1}{2}}\le\int_0^{\frac{1}{n}}(1-p^2)^{\frac{n^2-2n-1}{2}}\,\mathrm{d}p\le\frac{1}{n}$$

and we get

$$c_n/2\le\mathbb{P}\left\{\mathrm{d}\big(0,\mathrm{aff}(\xi_1,\dots,\xi_n)\big)\le 1/n\right\}\le c_n,\tag{51}$$

where

$$c_n=\frac{|S^{n-2}|^n}{n|S^{n-1}|^{n-1}}\frac{\Gamma(\frac{(n-1)^2+1}{2})}{\Gamma(\frac{(n-1)^2}{2})}\left(\frac{\Gamma(\frac{n-1}{2})}{\Gamma(\frac{n}{2})}\right)^{n-1}\frac{\Gamma(\frac{n-1}{2})}{\Gamma(\frac{1}{2})}.\tag{52}$$

Now we estimate

$$\mathbb{P}\left\{\mathrm{d}\big(\xi_{n+1},\mathrm{aff}(\xi_1,\dots,\xi_n)\big)\le\frac{nt}{\mathrm{vol}_{n-1}([\xi_1,\dots,\xi_n])}\,\middle|\,\mathrm{d}\big(0,\mathrm{aff}(\xi_1,\dots,\xi_n)\big)\le\frac{1}{n}\right\}.\tag{53}$$

First we estimate this expression from below for the range $0<t\le\frac{1}{n\cdot n!}$. The simplex with greatest volume in a Euclidean ball is a regular simplex. The volume of a regular simplex in $\mathbb{R}^m$ with sidelength $s$ is $\frac{\sqrt{m+1}}{m!}\left(\frac{s}{\sqrt{2}}\right)^m$ and the distance of its centroid to a vertex is $\sqrt{\frac{m}{2(m+1)}}s$. This means that the largest regular simplex inside a Euclidean ball with radius 1 has sidelength $\sqrt{2(m+1)/m}$ and volume $\frac{\sqrt{m+1}}{m!}\left(\frac{m+1}{m}\right)^{m/2}$. Therefore,

$$\mathrm{vol}_{n-1}([\xi_1,\dots,\xi_n])\le\frac{\sqrt{n}}{(n-1)!}\left(\frac{n}{n-1}\right)^{\frac{n-1}{2}}$$

and, for $M:=\frac{tn!}{\sqrt{n}}\left(\frac{n-1}{n}\right)^{\frac{n-1}{2}}$, (53) is greater than

$$\mathbb{P}\left\{\mathrm{d}\big(\xi_{n+1},\mathrm{aff}(\xi_1,\dots,\xi_n)\big)\le M\,\big|\,\mathrm{d}\big(0,\mathrm{aff}(\xi_1,\dots,\xi_n)\big)\le 1/n\right\}.\tag{54}$$

For each choice $\xi_1,\dots,\xi_n$ there is a hyperplane $H(\theta,p)$ spanned by this choice, $p$ is its distance from the origin and $\theta$ the normal to $H(\theta,p)$. Then we may choose $\xi_{n+1}$ from the set

$$S^{n-1}\cap H^+(\theta,p+M)\cap H^-(\theta,p-M),\tag{55}$$

where $H^+(\theta, p+M)$ and $H^-(\theta, p-M)$ are the halfspaces containing $p\,\theta$. By the condition $\mathrm{d}\big(0, \mathrm{aff}(\xi_1,\dots,\xi_n)\big) \le \frac{1}{n}$ we have $|p| \le \frac{1}{n}$ and by the assumption $t \le \frac{1}{n\cdot n!}$ and as $n \ge 3$ we get

$$p + M \le \frac{1}{n} + \frac{1}{n^{3/2}} \le \frac{2}{3}. \tag{56}$$

The surface area of (55) is greater than that of the right cylinder (without top and bottom) $\sqrt{1-(p+M)^2}S^{n-2} \times [p-M, p+M]$. The surface area of the cylinder equals $|S^{n-2}|\left(1-(p+M)^2\right)^{\frac{n-2}{2}} 2M$. By (56), we can apply Bernoulli's inequality and the surface area is at least

$$|S^{n-2}|\left(1 - \frac{n-2}{2}(p+M)^2\right)2M \ge \frac{|S^{n-2}|}{\sqrt{e}}\left(1 - \frac{n-2}{2}\left(\frac{1}{n} + \frac{tn!}{\sqrt{n}}\right)^2\right)\frac{2tn!}{\sqrt{n}},$$

where we used $\left(\frac{n}{n-1}\right)^{n-1} \le e$. Therefore, for all $t$ with $0 \le t \le \frac{1}{n\cdot n!}$ we may choose $\xi_{n+1}$ from a set of $(n-1)$-dimensional volume greater than

$$\frac{|S^{n-2}|}{\sqrt{e}}\frac{tn!}{\sqrt{n}}.$$

It follows, for all $t$ with $0 \le t \le \frac{1}{n\cdot n!}$

$$\begin{aligned}&\mathbb{P}\left\{\mathrm{d}\big(\xi_{n+1}, \mathrm{aff}(\xi_1,\dots,\xi_n)\big) \le \frac{nt}{\mathrm{vol}_{n-1}([\xi_1,\dots,\xi_n])}\,\middle|\, \mathrm{d}\big(0, \mathrm{aff}(\xi_1,\dots,\xi_n)\big) \le \frac{1}{n}\right\}\\ &\ge \frac{tn!}{\sqrt{e}\sqrt{n}}\frac{|S^{n-2}|}{|S^{n-1}|}.\end{aligned} \tag{57}$$

By (49), (51), (52) and (57),

$$\mathbb{P}\big(\mathrm{vol}_n([\xi_1,\dots,\xi_{n+1}]) \le t\big) \ge t\left(\frac{n!}{\sqrt{e}\sqrt{n}}\frac{|S^{n-2}|^{n+1}}{n|S^{n-1}|^n}\frac{\Gamma(\frac{(n-1)^2+1}{2})}{\Gamma(\frac{(n-1)^2}{2})}\left(\frac{\Gamma(\frac{n-1}{2})}{\Gamma(\frac{n}{2})}\right)^{n-1}\frac{\Gamma(\frac{n-1}{2})}{\Gamma(\frac{1}{2})}\right).$$

We have established the left hand side of inequality (37) for the range $0 \le t \le \frac{1}{n\cdot n!}$. We extend the range to $0 \le t \le |B_2^n|$. Let $\frac{1}{n\cdot n!} \le t \le |B_2^n|$. We have

$$\begin{aligned}\mathbb{P}\big(\mathrm{vol}_n([\xi_1,\dots,\xi_{n+1}]) \le t\big) &\ge \mathbb{P}\left(\mathrm{vol}_n([\xi_1,\dots,\xi_{n+1}]) \le \frac{1}{nn!}\right)\\ &\ge \frac{1}{nn!}\left(\frac{n!}{\sqrt{e}\sqrt{n}}\frac{|S^{n-2}|^{n+1}}{n|S^{n-1}|^n}\frac{\Gamma(\frac{(n-1)^2+1}{2})}{\Gamma(\frac{(n-1)^2}{2})}\left(\frac{\Gamma(\frac{n-1}{2})}{\Gamma(\frac{n}{2})}\right)^{n-1}\frac{\Gamma(\frac{n-1}{2})}{\Gamma(\frac{1}{2})}\right)\\ &\ge \frac{t}{|B_2^n|}\left(\frac{1}{\sqrt{e}\sqrt{n}}\frac{|S^{n-2}|^{n+1}}{n^2|S^{n-1}|^n}\frac{\Gamma(\frac{(n-1)^2+1}{2})}{\Gamma(\frac{(n-1)^2}{2})}\left(\frac{\Gamma(\frac{n-1}{2})}{\Gamma(\frac{n}{2})}\right)^{n-1}\frac{\Gamma(\frac{n-1}{2})}{\Gamma(\frac{1}{2})}\right).\end{aligned}$$

Now we show the right hand side estimate of (37). We have

$$\begin{aligned}&\mathbb{P}\big(\mathrm{vol}_n([\xi_1,\dots,\xi_{n+1}]) \le t\big) = \mathbb{P}\left\{\mathrm{d}\big(\xi_{n+1}, \mathrm{aff}(\xi_1,\dots,\xi_n)\big) \le \frac{nt}{\mathrm{vol}_{n-1}([\xi_1,\dots,\xi_n])}\right\}\\ &= \frac{1}{|S^{n-1}|^n}\int_{(S^{n-1})^n}\mathbb{P}\left\{\xi_{n+1} : \mathrm{d}\big(\xi_{n+1}, \mathrm{aff}(\xi_1,\dots,\xi_n)\big) \le \frac{nt}{\mathrm{vol}_{n-1}([\xi_1,\dots,\xi_n])}\right\}\mathrm{d}\xi_1\cdots\mathrm{d}\xi_n.\end{aligned}$$

We estimate the integrand. Let $H(\theta,p)$ be the hyperplane spanned by $\xi_1,\dots,\xi_n$ where $\theta$ is its unit normal and $p$ its distance to 0. The set of all $\xi_{n+1}$ with $\mathrm{d}\big(\xi_{n+1},\mathrm{aff}(\xi_1,\dots,\xi_n)\big)\leq \frac{nt}{\mathrm{vol}_{n-1}([\xi_1,\dots,\xi_n])}$ is

$$S^{n-1}\cap H^+\left(\theta, p+\frac{nt}{\mathrm{vol}_{n-1}([\xi_1,\dots,\xi_n])}\right)\cap H^-\left(\theta, p-\frac{nt}{\mathrm{vol}_{n-1}([\xi_1,\dots,\xi_n])}\right),$$

where both half spaces contain the point $p\,\theta$. By Lemma 3 the above set has area

$$|S^{n-2}|\int_{\max\{-1,p-\frac{nt}{\mathrm{vol}_{n-1}([\xi_1,\dots,\xi_n])}\}}^{\min\{1,p+\frac{nt}{\mathrm{vol}_{n-1}([\xi_1,\dots,\xi_n])}\}}(1-s^2)^{\frac{n-3}{2}}\,\mathrm{d}s\leq 2|S^{n-2}|\frac{nt}{\mathrm{vol}_{n-1}([\xi_1,\dots,\xi_n])}.$$

where we used that the integrand is smaller than 1. Thus,

$$\mathbb{P}\left\{\xi_{n+1}:\mathrm{d}\big(\xi_{n+1},\mathrm{aff}(\xi_1,\dots,\xi_n)\big)\leq\frac{nt}{\mathrm{vol}_{n-1}([\xi_1,\dots,\xi_n])}\right\}\leq\frac{|S^{n-2}|}{|S^{n-1}|}\frac{2nt}{\mathrm{vol}_{n-1}([\xi_1,\dots,\xi_n])}.$$

Altogether, using the formula of Blaschke-Petkantschin (9),

$$\begin{aligned}\mathbb{P}\big(\mathrm{vol}_n([\xi_1,\dots,\xi_{n+1}])\leq t\big)&\leq\frac{2nt|S^{n-2}|}{|S^{n-1}|^{n+1}}\int_{S^{n-1}}\cdots\int_{S^{n-1}}\frac{\mathrm{d}\xi_1\cdots\mathrm{d}\xi_n}{\mathrm{vol}_{n-1}([\xi_1,\dots,\xi_n])}\\&=\frac{2n!t|S^{n-2}|}{|S^{n-1}|^{n+1}}\int_0^1\int_{S^{n-1}}\int_{(H(\theta,p)\cap S^{n-1})^n}\frac{1}{(1-p^2)^{n/2}}\,\mathrm{d}\xi_1\cdots\mathrm{d}\xi_n\,\mathrm{d}\theta\,\mathrm{d}p\\&=\frac{2n!t|S^{n-2}|}{|S^{n-1}|^{n}}\int_0^1\int_{(H(\theta,p)\cap S^{n-1})^n}\frac{1}{(1-p^2)^{n/2}}\,\mathrm{d}\xi_1\cdots\mathrm{d}\xi_n\,\mathrm{d}p.\end{aligned}$$

The set $H(\theta,p)\cap S^{n-1}$ is an $(n-2)$-dimensional Euclidean sphere with radius $\sqrt{1-p^2}$. Therefore, the previous expression equals

$$\begin{aligned}\frac{2n!t|S^{n-2}|^{n+1}}{|S^{n-1}|^n}\int_0^1(1-p^2)^{\frac{n^2-3n}{2}}\,\mathrm{d}p&=\frac{2n!t|S^{n-2}|^{n+1}}{|S^{n-1}|^n}\frac{1}{2}\int_0^1(1-w)^{\frac{n^2-3n}{2}}w^{-\frac12}\,\mathrm{d}p\\&=\frac{n!t|S^{n-2}|^{n+1}}{|S^{n-1}|^n}B\left(\frac12,\frac{n^2-3n+2}{2}\right).\end{aligned}$$

The claim follows. □

We define a function that is instrumental to estimate the expected minimal volume facet. For the estimate of the expected maximal volume facet we introduced a similar function (18). Its difference reflects that we consider the minimal facet and the maximal facet. Let $t\geq 0$ and let $\tilde{H}_{t,C}:(S^{n-1})^n\to\{0,1\}$ be such that $\tilde{H}_{t,C}(\xi_1,\dots,\xi_n)=1$ if the following holds:

$$\begin{aligned}&\mathrm{aff}(\xi_1,\dots,\xi_n)\cap S^{n-1}\subseteq C\quad\text{and}\\&\mathrm{vol}_{n-1}([\xi_1,\dots,\xi_n])\leq t\left(1-d^2(0,\mathrm{aff}(\xi_1,\dots,\xi_n))\right)^{\frac{n-1}{2}}.\end{aligned}\tag{58}$$

Else we put $\tilde{H}_{t,C}(\xi_1,\dots,\xi_n)=0$.

**Lemma 15.** *Let $C$ be a cap of $S^{n-1}$ with $\mathrm{vol}_{n-1}(C) = \frac{|S^{n-1}|}{N}$. Then*
*(i) For $n = 3$,*

$$\frac{4}{\pi^3 10^5}\frac{t^{5/3}}{N^3} \le \mathbb{E}\,\tilde{H}_{t,C} \le 171\frac{\pi^7}{16}\frac{t^{5/3}}{N^3}, \tag{59}$$

*where the left hand side inequality holds for all $t$ with $0 \le t \le \pi$ and the right hand side inequality for all $t$ with $0 \le t$.*
*(ii) For all $n \ge 4$,*

$$\frac{t^2}{N^n}\frac{3a_{n-1}^2(n-1)!(n-1)^{n-3}|S^{n-2}|}{32\cdot b_{n-1}\pi^{n^2-n-2}} \le \mathbb{E}\,\tilde{H}_{t,C} \le \frac{t^2}{N^n}b_{n-1}\left(\frac{\pi}{2}\right)^{n^2-n}(n-1)^{n-3}(n-1)!|S^{n-2}|, \tag{60}$$

*where $0 \le t \le |B_2^{n-1}|$ for the left hand side inequality and $0 \le t$ for the right hand side inequality. The constant $a_{n-1}$ is given by* (38) *and $b_{n-1}$ by* (39)*.*

*Proof.* We proceed as in the proof of Lemma 13. By definition

$$\mathbb{E}\,\tilde{H}_{t,C} = \frac{1}{|S^{n-1}|^n}\int_{S^{n-1}}\cdots\int_{S^{n-1}}\tilde{H}_{t,C}(\xi_1,\dots,\xi_n)\,\mathrm{d}\xi_1\cdots\mathrm{d}\xi_n.$$

We apply the spherical Blaschke-Petkantschin formula (9),

$$\mathbb{E}\,\tilde{H}_{t,C} = \frac{(n-1)!}{|S^{n-1}|^n}\int_0^1\int_{S^{n-1}}\int_{(H(\theta,p)\cap S^{n-1})^n}\tilde{H}_{t,C}(\xi_1,\dots,\xi_n)\frac{\big|[\xi_1,\dots,\xi_n]\big|}{(1-p^2)^{n/2}}\,\mathrm{d}\xi_1\cdots\mathrm{d}\xi_n\,\mathrm{d}\theta\,\mathrm{d}p. \tag{61}$$

Let $\phi(N,n)$ be the angle of the cap $C$. Choosing $R = N$ in Lemma 4,

$$L := \left(\frac{|S^{n-1}|}{N|B_2^{n-1}|}\right)^{\frac{1}{n-1}} \le \phi(N,n) \le U := \frac{\pi}{2}\left(\frac{|S^{n-1}|}{N|B_2^{n-1}|}\right)^{\frac{1}{n-1}}. \tag{62}$$

As in the proof of Lemma 13, the first condition in the definition of $\tilde{H}_{t,C}$ translates into $\arccos\langle\theta,e_1\rangle + \arccos q \le \phi(N,n)$ (see (20)) and we get

$$\mathbb{E}\,\tilde{H}_{t,C} = \frac{(n-1)!|S^{n-2}|^n}{|S^{n-1}|^n}\int_0^1\int_{S^{n-1}}\chi\big(\arccos\langle\theta,e_1\rangle + \arccos q \le \phi(N,n)\big) \\ I(q,t)(1-q^2)^{\frac{n^2-2n-1}{2}}\,\mathrm{d}\theta\,\mathrm{d}q, \tag{63}$$

where we assume that the hyperplane $H$ with $C = S^{n-1}\cap H^-$ is orthogonal to $e_1$. Moreover, using the substitution $\xi_i = (1-q^2)^{\frac{1}{2}}\zeta_i$,

$$\begin{aligned} I(q,t) &= \int_{(H(e_1,q)\cap S^{n-1})^n}\chi\left(\mathrm{vol}_{n-1}([\xi_1,\dots,\xi_n]) \le t(1-q^2)^{\frac{n-1}{2}}\right)\frac{\big|[\xi_1,\dots,\xi_n]\big|}{(1-q^2)^{\frac{n-1}{2}}}\prod_{i=1}^n\frac{\mathrm{d}\xi_i}{|S^{n-2}|(1-q^2)^{\frac{n-2}{2}}} \\ &= \frac{1}{|S^{n-2}|^n}\int_{(S^{n-2})^n}\chi\big(\mathrm{vol}_{n-1}([\zeta_1,\dots,\zeta_n]) \le t\big)\,\mathrm{vol}_{n-1}([\zeta_1,\dots,\zeta_n])\,\mathrm{d}\zeta_1\cdots\mathrm{d}\zeta_n \\ &= \mathbb{E}\Big(\chi\big(\mathrm{vol}_{n-1}([\zeta_1,\dots,\zeta_n]) \le t\big)\,\mathrm{vol}_{n-1}([\zeta_1,\dots,\zeta_n])\Big). \end{aligned}$$

Thus

$$
\begin{aligned}
I(q,t) &= \int_0^\infty \mathbb{P}\Big(\chi\big(\mathrm{vol}_{n-1}([\zeta_1,\dots,\zeta_n]) \le t\big)\,\mathrm{vol}_{n-1}([\zeta_1,\dots,\zeta_n]) \ge s\Big)\,\mathrm{d}s \\
&= \int_0^t \mathbb{P}\big(s \le \mathrm{vol}_{n-1}([\zeta_1,\dots,\zeta_n]) \le t\big)\,\mathrm{d}s. \qquad (64)
\end{aligned}
$$

**Upper bound on $\mathbb{E}\,\tilde{H}_{t,C}$.** In order to obtain an upper bound on $\mathbb{E}\,\tilde{H}_{t,C}$ we first upper bound $I(q,t)$. By (64), and Lemma 14 applied to the dimension $n-1 \ge 2$, we get for all $t \ge 0$,

$$
I(q,t) \le t \cdot \mathbb{P}\big(\mathrm{vol}_{n-1}([\zeta_1,\dots,\zeta_n]) \le t\big) \le \bar{I}_{n-1}(t), \qquad (65)
$$

where $\bar{I}_n(t) = b_n\,t^2$ if $n \ge 3$, $b_n$ is given by (39), and $\bar{I}_2(t) = 342\,t^{\frac{5}{3}}$ by (36). By (63) and (65), and then by (22),

$$
\begin{aligned}
\mathbb{E}\,\tilde{H}_{t,C} &\le \bar{I}_{n-1}(t)\frac{(n-1)!|S^{n-2}|^n}{|S^{n-1}|^n} \\
&\qquad \int_0^1 \int_{S^{n-1}} \chi\big(\arccos\langle\theta,e_1\rangle + \arccos q \le \phi(N,n)\big)(1-q^2)^{\frac{n^2-2n-1}{2}}\,\mathrm{d}\theta\,\mathrm{d}q \\
&\le \bar{I}_{n-1}(t)\frac{(n-1)!|S^{n-2}|^{n+1}}{|S^{n-1}|^n} \\
&\qquad \int_0^1 \int_0^\pi \chi\big(\phi_1 + \arccos q \le \phi(N,n)\big)(1-q^2)^{\frac{n^2-2n-1}{2}} \sin^{n-2}\phi_1\,\mathrm{d}\phi_1\,\mathrm{d}q.
\end{aligned} \qquad (66)
$$

By (62) the double integral in (66) is at most

$$
\begin{aligned}
&\int_0^1 \int_0^\pi \chi\left(\phi_1 + \arccos q \le U\right)(1-q^2)^{\frac{n^2-2n-1}{2}} \sin^{n-2}\phi_1\,\mathrm{d}\phi_1\,\mathrm{d}q \\
&\le \int_0^1 \chi\left(\arccos q \le U\right)(1-q^2)^{\frac{n^2-2n-1}{2}}\,\mathrm{d}q \int_0^\pi \chi\left(\phi_1 \le U\right)\sin^{n-2}\phi_1\,\mathrm{d}\phi_1.
\end{aligned}
$$

Therefore

$$
\mathbb{E}\,\tilde{H}_{t,C} \le \bar{I}_{n-1}(t)\frac{(n-1)!|S^{n-2}|^{n+1}}{|S^{n-1}|^n}\int_{\cos(U)}^1 (1-q^2)^{\frac{n^2-2n-1}{2}}\,\mathrm{d}q \int_0^U \sin^{n-2}\phi_1\,\mathrm{d}\phi_1. \qquad (67)
$$

Since $\sin\phi_1 \le \phi_1$, we have $\int_0^U \sin^{n-2}\phi_1\,\mathrm{d}\phi_1 \le \int_0^U \phi_1^{n-2}\,\mathrm{d}\phi_1 = U^{n-1}/(n-1)$. With $q = \cos\alpha$ and $\mathrm{d}q = -\sin\alpha\,\mathrm{d}\alpha$,

$$
\int_{\cos(U)}^1 (1-q^2)^{\frac{n^2-2n-1}{2}}\,\mathrm{d}q = \int_0^U (\sin\alpha)^{n^2-2n} d\alpha \le \int_0^U \alpha^{n^2-2n} d\alpha = \frac{U^{(n-1)^2}}{(n-1)^2}
$$

Altogether, using the previous estimates in (67), we get for $n \ge 3$ and $t \ge 0$,

$$
\begin{aligned}
\mathbb{E}\,\tilde{H}_{t,C} &\le \bar{I}_{n-1}(t)\frac{(n-1)!|S^{n-2}|^{n+1}}{|S^{n-1}|^n}\frac{1}{N}\left(\frac{\pi}{2}\right)^{n-1}\frac{|S^{n-1}|}{|S^{n-2}|}N^{-n+1}\frac{\left(\frac{\pi}{2}\right)^{(n-1)^2}}{(n-1)^2}\left(\frac{|S^{n-1}|}{|B_2^{n-1}|}\right)^{n-1} \\
&\le \bar{I}_{n-1}(t)N^{-n}|S^{n-2}|\left(\pi/2\right)^{n^2-n}(n-1)^{n-3}(n-1)!\ .
\end{aligned}
$$

The claimed upper bounds (Equations (59) and (60)) follow.

**Lower bound on $\mathbb{E}\,\tilde{H}_{t,C}$.** Similarly to the upper bound, first we get a lower bound on $I(q,t)$. For $n \geq 4$ and by (64),

$$I(q,t) = \int_0^t \mathbb{P}\big(\mathrm{vol}_{n-1}([\zeta_1,\ldots,\zeta_n]) \leq t\big) - \mathbb{P}\big(\mathrm{vol}_{n-1}([\zeta_1,\ldots,\zeta_n]) \leq s\big)\,\mathrm{d}s.$$

By (37) we have $a_{n-1} \leq b_{n-1}$. Therefore $\frac{a_{n-1}}{2b_{n-1}}t \leq t$ and

$$I(q,t) \geq \int_0^{\frac{a_{n-1}}{2b_{n-1}}t} \mathbb{P}\big(\mathrm{vol}_{n-1}([\zeta_1,\ldots,\zeta_n]) \leq t\big) - \mathbb{P}\big(\mathrm{vol}_{n-1}([\zeta_1,\ldots,\zeta_n]) \leq s\big)\,\mathrm{d}s.$$

We apply (37). For $t$ with $0 \leq t \leq |B_2^{n-1}|$,

$$I(q,t) \geq \int_0^{\frac{a_{n-1}}{2b_{n-1}}t} a_{n-1}t - b_{n-1}s\,\mathrm{d}s = a_{n-1}ts - \frac{1}{2}b_{n-1}s^2\Big|_0^{\frac{a_{n-1}}{2b_{n-1}}t} = \frac{3a_{n-1}^2}{8b_{n-1}}t^2 := \underline{I}_{n-1}(t).$$

Now we get a lower bound on $I(q,t)$ for $n=3$. By (64) and with $c := 9/10^5$,

$$I(q,t) \geq \int_0^{ct} \mathbb{P}\big(\mathrm{vol}_{n-1}([\zeta_1,\zeta_2,\zeta_3]) \leq t\big) - \mathbb{P}\big(\mathrm{vol}_{n-1}([\zeta_1,\zeta_2,\zeta_3]) \leq s\big)\,\mathrm{d}s.$$

By (36),

$$I(q,t) \geq \int_0^{ct} \frac{(2t)^{\frac{2}{3}}}{\pi^{\frac{8}{3}}} - 342 s^{\frac{2}{3}}\,\mathrm{d}s \geq ct^{\frac{5}{3}}\left(\frac{2^{\frac{2}{3}}}{\pi^{\frac{8}{3}}} - 206\cdot c^{\frac{2}{3}}\right) \geq \frac{4}{10^5}t^{\frac{5}{3}} := \underline{I}_2(t).$$

It follows from this and (63) that for $t$ with $0 \leq t \leq |B_2^{n-1}|$,

$$\mathbb{E}\,\tilde{H}_{t,C} \geq \underline{I}_{n-1}(t)\frac{(n-1)!|S^{n-2}|^n}{|S^{n-1}|^n}$$
$$\int_0^1\int_{S^{n-1}} \chi\big(\arccos\langle\theta,e_1\rangle + \arccos q \leq \phi(N,n)\big)(1-q^2)^{\frac{n^2-2n-1}{2}}\,\mathrm{d}\theta\,\mathrm{d}q.$$

By (22), for $t$ with $0 \leq t \leq |B_2^{n-1}|$

$$\mathbb{E}\,\tilde{H}_{t,C} \geq \underline{I}_{n-1}(t)\frac{(n-1)!|S^{n-2}|^{n+1}}{|S^{n-1}|^n}$$
$$\int_0^1\int_0^{\pi} \chi\big(\phi_1 + \arccos q \leq \phi(N,n)\big)(1-q^2)^{\frac{n^2-2n-1}{2}}\sin^{n-2}\phi_1\,\mathrm{d}\phi_1\,\mathrm{d}q.$$

By (62) the double integral is greater than

$$\int_0^1\int_0^{\pi} \chi\left(\phi_1 + \arccos q \leq L\right)(1-q^2)^{\frac{n^2-2n-1}{2}}\sin^{n-2}\phi_1\,\mathrm{d}\phi_1\,\mathrm{d}q.$$

Since $\chi\left(\phi_1 + \arccos q \leq u\right) \geq \chi\left(\phi_1 \leq u/2\right)\chi\left(\arccos q \leq u/2\right)$, the double integral is greater than

$$\int_0^1 \chi\left(\arccos q \leq \frac{L}{2}\right)(1-q^2)^{\frac{n^2-2n-1}{2}}\,\mathrm{d}q\int_0^{\pi}\chi\left(\phi_1 \leq \frac{L}{2}\right)\sin^{n-2}\phi_1\,\mathrm{d}\phi_1.$$

We have $L/2 \leq \pi$. Therefore, for $t$ with $0 \leq t \leq |B_2^n|$,

$$\mathbb{E}\,\tilde{H}_{t,C} \geq \underline{I}_{n-1}(t) \frac{(n-1)!|S^{n-2}|^{n+1}}{|S^{n-1}|^n} \int_{\cos(L/2)}^1 (1-q^2)^{\frac{n^2-2n-1}{2}} \,\mathrm{d}q \int_0^{L/2} \sin^{n-2}\phi_1 \,\mathrm{d}\phi_1. \tag{68}$$

Since $\sin\phi_1 \geq \frac{2}{\pi}\phi_1$,

$$\int_0^{L/2} \sin^{n-2}\phi_1 \,\mathrm{d}\phi_1 \geq \left(\frac{2}{\pi}\right)^{n-2} \int_0^{L/2} \phi_1^{n-2} \,\mathrm{d}\phi_1 = \frac{1}{2N(n-1)} \left(\frac{1}{\pi}\right)^{n-2} \frac{|S^{n-1}|}{|B_2^{n-1}|}. \tag{69}$$

With $q = \cos\alpha$ and $\mathrm{d}q = -\sin\alpha \,\mathrm{d}\alpha$, $\int_{\cos(L/2)}^1 (1-q^2)^{\frac{n^2-2n-1}{2}} \,\mathrm{d}q = \int_0^{L/2} (\sin\alpha)^{n^2-2n} \,\mathrm{d}\alpha$. Using again $\sin\alpha \geq \frac{2}{\pi}\alpha$ for $0 \leq \alpha \leq \pi/2$, this is greater than

$$\begin{aligned}\left(\frac{2}{\pi}\right)^{n^2-2n} \int_0^{\frac{L}{2}} \alpha^{n^2-2n} d\alpha &= \frac{1}{(n-1)^2} \left(\frac{2}{\pi}\right)^{n^2-2n} \left(\frac{L}{2}\right)^{(n-1)^2} \\ &= N^{-n+1} \left(\frac{1}{\pi}\right)^{n^2-2n} \frac{1}{2(n-1)^2} \left(\frac{|S^{n-1}|}{|B_2^{n-1}|}\right)^{n-1}.\end{aligned} \tag{70}$$

Altogether, by (68), (69) and (70),

$$\begin{aligned}\mathbb{E}\,\tilde{H}_{t,C} &\geq \underline{I}_{n-1}(t) \frac{(n-1)!|S^{n-2}|^{n+1}}{|S^{n-1}|^n} \frac{1}{2N} \left(\frac{1}{\pi}\right)^{n^2-n-2} \frac{|S^{n-1}|}{|S^{n-2}|} N^{-n+1} \frac{1}{2(n-1)^2} \left(\frac{|S^{n-1}|}{|B_2^{n-1}|}\right)^{n-1} \\ &= \frac{\underline{I}_{n-1}(t)}{N^n} \frac{(n-1)!(n-1)^{n-3}|S^{n-2}|}{4\pi^{n^2-n-2}}.\end{aligned}$$

□

# 7 Expected minimum upper bound, $n \geq 3$

**Proof of the upper estimates in Theorem 2, $n \geq 3$.** Let $N \in \mathbb{N}$ and let $C_1, \dots, C_k$ be a maximal set of caps in $S^{n-1}$ such that

$$\mathrm{vol}_{n-1}(C_i) = \mathrm{vol}_{n-1}(S^{n-1})/N \qquad i = 1, \dots, k \tag{71}$$

and such that for all $i \neq j$, $\mathrm{int}(C_i) \cap \mathrm{int}(C_j) = \emptyset$. By Lemma 5, $3^{-n}N \leq k \leq N$. We are now defining the function $\tilde{G}_{i,t} : S^{n-1} \times \cdots \times S^{n-1} \to \{0,1\}$. The definition of $\tilde{G}_{i,t}$ is very similar to the definition of $G_i$. Their difference lies in conditions (74) and (26) . Let $t \geq 0$. We define $\tilde{G}_{i,t}(\xi_1, \dots, \xi_N) = 1$ if the following conditions hold:

$$\text{Exactly } n \text{ points } \xi_{\ell_1}, \dots, \xi_{\ell_n} \text{ are chosen from } C_i \text{ and they are affinely independent,} \tag{72}$$

$$\mathrm{aff}(\xi_{\ell_1}, \dots, \xi_{\ell_n}) \cap S^{n-1} \subseteq C_i. \tag{73}$$

$$\mathrm{vol}_{n-1}([\xi_{\ell_1}, \dots, \xi_{\ell_n}]) \leq t\big(1 - d^2(0, \mathrm{aff}(\xi_{\ell_1}, \dots, \xi_{\ell_n}))\big)^{\frac{n-1}{2}}. \tag{74}$$

Else we put $\tilde{G}_{i,t}(\xi_1, \dots, \xi_N) = 0$.

In the same way as we established (32) (with $R = N$) we arrive now at

$$\mathbb{P}\left(\sum_{i=1}^{k} \tilde{G}_{i,t} > 0\right) \geq \frac{\left[k\binom{N}{n}\left(1-\frac{1}{N}\right)^{N-n}\mathbb{E}\,\tilde{H}_{t,C_1}\right]^2}{k(k-1)\binom{N}{n}\binom{N-n}{n}\left(1-\frac{2}{N}\right)^{N-2n}(\mathbb{E}\,\tilde{H}_{t,C_1})^2 + k\binom{N}{n}\left(1-\frac{1}{N}\right)^{N-n}\mathbb{E}\,\tilde{H}_{t,C_1}}$$

$$= \frac{k\binom{N}{n}\left(1-\frac{1}{N}\right)^{2N-2n}\mathbb{E}\,\tilde{H}_{t,C_1}}{(k-1)\binom{N-n}{n}\left(1-\frac{2}{N}\right)^{N-2n}\mathbb{E}\,\tilde{H}_{t,C_1} + \left(1-\frac{1}{N}\right)^{N-n}}. \tag{75}$$

We will upper bound the reciprocal of this. For $N \geq 4n+1$ we have

$$\frac{\left(1-\frac{2}{N}\right)^{N-2n}}{\left(1-\frac{1}{N}\right)^{2(N-n)}} = \left(\frac{1-\frac{2}{N}}{1-\frac{1}{N}}\right)^{N-2n}\left(1-\frac{1}{N}\right)^{-N} = \left(1-\frac{1}{N-1}\right)^{N-2n}\left(1-\frac{1}{N}\right)^{-N}$$

$$= \left(1-\frac{1}{(N-1)^2}\right)^{N}\left(1-\frac{1}{N-1}\right)^{-2n} \leq \left(1-\frac{2n}{N-1}\right)^{-1} \leq 1+\frac{4n}{N-1}.$$

Thus, for $N \geq 4n+1$ and using $\binom{N-n}{n}/\binom{N}{n} \leq 1$ we have

$$\frac{(k-1)\binom{N-n}{n}\left(1-\frac{2}{N}\right)^{N-2n}\mathbb{E}\,\tilde{H}_{t,C_1}}{k\binom{N}{n}\left(1-\frac{1}{N}\right)^{2(N-n)}\mathbb{E}\,\tilde{H}_{t,C_1}} \leq 1+\frac{4n}{N-1}. \tag{76}$$

Moreover, since $(1-\frac{1}{N})^{N-n} \geq (1-\frac{1}{N})^{N-1} = \frac{1}{(1+\frac{1}{N-1})^{N-1}} > 1/e$,

$$\frac{\left(1-\frac{1}{N}\right)^{N-n}}{k\binom{N}{n}\left(1-\frac{1}{N}\right)^{2(N-n)}\mathbb{E}\,\tilde{H}_{t,C_1}} = \frac{1}{k\binom{N}{n}\left(1-\frac{1}{N}\right)^{N-n}\mathbb{E}\,\tilde{H}_{t,C_1}} \leq \frac{e}{k\binom{N}{n}\mathbb{E}\,\tilde{H}_{t,C_1}}. \tag{77}$$

By (7) we have $k \geq 3^{-n}N$. This and the estimate (60) for $\mathbb{E}\,\tilde{H}_{t,C_1}$ gives for $n \geq 4$,

$$\frac{\left(1-\frac{1}{N}\right)^{N-n}}{k\binom{N}{n}\left(1-\frac{1}{N}\right)^{2(N-n)}\mathbb{E}\,\tilde{H}_{t,C_1}} \leq \frac{e\cdot 3^n}{N\binom{N}{n}}\frac{N^n}{t^2}\frac{32\cdot b_{n-1}\pi^{n^2-n-2}}{3a_{n-1}^2(n-1)!(n-1)^{n-3}|S^{n-2}|}$$

$$\leq \frac{1}{t^2N}\frac{32e\cdot 3^{n-1}ne^{\frac{n^2}{N-n}}\cdot b_{n-1}\pi^{n^2-n-2}}{a_{n-1}^2(n-1)^{n-3}|S^{n-2}|}, \tag{78}$$

where $a_{n-1}$ and $b_{n-1}$ are defined in (38) and (39). Similarly, by (59) for $n = 3$,

$$\frac{\left(1-\frac{1}{N}\right)^{N-3}}{k\binom{N}{3}\left(1-\frac{1}{N}\right)^{2(N-3)}\mathbb{E}\,\tilde{H}_{t,C_1}} \leq 10^5\frac{N^3}{t^{5/3}}\frac{e(3\pi)^3}{4N\binom{N}{3}} \leq \frac{30^5\pi^3e^{\frac{9}{N-3}}}{t^{5/3}N}. \tag{79}$$

By (75), (76) and (78) for dimensions greater than or equal to 4 and sufficiently large $N$,

$$\mathbb{P}\left(\sum_{i=1}^{k} \tilde{G}_{i,t} > 0\right) \geq \left(1+\frac{4n}{N-1}+\frac{1}{t^2N}\frac{32e\cdot 3^{n-1}ne^{\frac{n^2}{N-n}}\cdot b_{n-1}\pi^{n^2-n-2}}{a_{n-1}^2(n-1)^{n-3}|S^{n-2}|}\right)^{-1} \geq 1-\frac{c_n}{Nt^2}, \tag{80}$$

where $c_n$ is a constant depending only on the dimension $n$. Similarly, by (75), (76) and (79) for dimension 3, we get with a suitably chosen constant $c_3$,

$$\mathbb{P}\left(\sum_{i=1}^{k}\tilde{G}_{i,t}>0\right)\geq\left(1+\frac{12}{N-1}+\frac{30^5\pi^3e^{\frac{9}{N-3}}}{t^{5/3}N}\right)^{-1}\geq 1-\frac{c_3}{Nt^{5/3}}. \tag{81}$$

We show now

$$\left\{(\xi_1,\dots,\xi_N):\min_{F\in\mathcal{F}(P_N)}\mathrm{vol}_{n-1}(F)\geq s\right\}\subseteq\left\{(\xi_1,\dots,\xi_N):\sum_{i=1}^{k}\tilde{G}_{i,s\frac{N|B_2^{n-1}|}{2|S^{n-1}|}}=0\right\}. \tag{82}$$

If $\sum_{i=1}^{k}\tilde{G}_{i,t}(\xi_1,\dots,\xi_N)>0$, then there is $i_0$ with $\tilde{G}_{i_0,t}(\xi_1,\dots,\xi_N)>0$. Since $\tilde{G}_{i_0,t}$ takes only the values 0 and 1, we have $\tilde{G}_{i_0,t}(\xi_1,\dots,\xi_N)=1$. This implies that $C_{i_0}$ contains $n$ points $\xi_{\ell_1},\dots,\xi_{\ell_n}$ whose convex hull $[\xi_{\ell_1},\dots,\xi_{\ell_n}]$ is a facet of $[\xi_1,\dots,\xi_N]$. Let $H$ be the hyperplane with $C_{i_0}=H^-\cap S^{n-1}$ and let $L$ be the hyperplane containing $\xi_{\ell_1},\dots,\xi_{\ell_n}$. By assumption (73), $L\cap S^{n-1}\subseteq C_{i_0}$. Consequently, by this and assumption (71),

$$\mathrm{vol}_{n-1}(L\cap S^{n-1})\leq\mathrm{vol}_{n-1}(C_{i_0})=|S^{n-1}|/N.$$

Thus with Lemma 3,

$$\left(1-d^2(0,\mathrm{aff}(\xi_{\ell_1},\dots,\xi_{\ell_n}))\right)^{\frac{n-1}{2}}|B_2^{n-1}|\leq\mathrm{vol}_{n-1}(L^-\cap S^{n-1})\leq|S^{n-1}|/N.$$

With (74), this gives

$$\mathrm{vol}_{n-1}([\xi_{\ell_1},\dots,\xi_{\ell_n}])\leq t\left(1-d^2(0,\mathrm{aff}(\xi_{\ell_1},\dots,\xi_{\ell_n}))\right)^{\frac{n-1}{2}}\leq t\frac{|S^{n-1}|}{N|B_2^{n-1}|}.$$

Altogether, if $\sum_{i=1}^{k}\tilde{G}_{i,t}>0$, then $[\xi_1,\dots,\xi_N]$ has a facet $[\xi_{\ell_1},\dots,\xi_{\ell_n}]$ with

$$\mathrm{vol}_{n-1}([\xi_{\ell_1},\dots,\xi_{\ell_n}])\leq t\frac{|S^{n-1}|}{N|B_2^{n-1}|}.$$

If $\min_{F\in\mathcal{F}(P_N)}\mathrm{vol}_{n-1}(F)\geq s$ and $\sum_{i=1}^{k}\tilde{G}_{i,t}>0$, then

$$s\leq t\frac{|S^{n-1}|}{N|B_2^{n-1}|}.$$

For $t=s\frac{N|B_2^{n-1}|}{2|S^{n-1}|}$ this estimate cannot hold. So we have shown that $\min_{F\in\mathcal{F}(P_N)}\mathrm{vol}_{n-1}(F)\geq s$ implies $\sum_{i=1}^{k}\tilde{G}_{i,s\frac{N|B_2^{n-1}|}{2|S^{n-1}|}}(\xi_1,\dots,\xi_N)=0$, i.e. we have shown (82). Consequently, by (82),

$$\mathbb{P}\left(\min_{F\in\mathcal{F}(P_N)}\mathrm{vol}_{n-1}(F)\geq s\right)\leq\mathbb{P}\left(\sum_{i=1}^{k}\tilde{G}_{i,s\frac{N|B_2^{n-1}|}{2|S^{n-1}|}}=0\right)=1-\mathbb{P}\left(\sum_{i=1}^{k}\tilde{G}_{i,s\frac{N|B_2^{n-1}|}{2|S^{n-1}|}}>0\right).$$

By this and (80) in the case of dimension greater than or equal to 4

$$\mathbb{P}\left(\min_{F\in\mathcal{F}(P_N)}\mathrm{vol}_{n-1}(F)\geq s\right)\leq\frac{4c_n|S^{n-1}|^2}{N^3s^2|B_2^{n-1}|^2} \tag{83}$$

and by (81) in dimension 3

$$\mathbb{P}\left(\min_{F\in\mathcal{F}(P_N)}\mathrm{vol}_{n-1}(F)\geq s\right)\leq\frac{32\cdot c_3}{N^{8/3}s^{5/3}}. \tag{84}$$

Since the polytope is contained in $B_2^n$, every facet has a surface area that is smaller than $|B_2^{n-1}|$. So we have

$$\begin{aligned}&\mathbb{E}\left(\min_{F\in\mathcal{F}(P_N)}\mathrm{vol}_{n-1}(F)\right)=\int_0^{|B_2^{n-1}|}\mathbb{P}\left(\min_{F\in\mathcal{F}(P_N)}\mathrm{vol}_{n-1}(F)\geq s\right)\mathrm{d}s\\&=\int_0^{N^{-3/2}}\mathbb{P}\left(\min_{F\in\mathcal{F}(P_N)}\mathrm{vol}_{n-1}(F)\geq s\right)\mathrm{d}s+\int_{N^{-3/2}}^{|B_2^{n-1}|}\mathbb{P}\left(\min_{F\in\mathcal{F}(P_N)}\mathrm{vol}_{n-1}(F)\geq s\right)\mathrm{d}s.\end{aligned}$$

By (83) with a new constant $c_n$ we get for $n\geq 4$

$$\mathbb{E}\left(\min_{F\in\mathcal{F}(P_N)}\mathrm{vol}_{n-1}(F)\right)\leq N^{-3/2}+\frac{c_n}{N^3}\int_{N^{-3/2}}^{|B_2^{n-1}|}\frac{1}{s^2}\,\mathrm{d}s\leq(1+c_n)N^{-3/2}.$$

In the case of dimension 3 we have

$$\begin{aligned}&\mathbb{E}\left(\min_{F\in\mathcal{F}(P_N)}\mathrm{vol}_{n-1}(F)\right)=\int_0^{\pi}\mathbb{P}\left(\min_{F\in\mathcal{F}(P_N)}\mathrm{vol}_{n-1}(F)\geq s\right)\mathrm{d}s\\&=\int_0^{N^{-8/5}}\mathbb{P}\left(\min_{F\in\mathcal{F}(P_N)}\mathrm{vol}_{n-1}(F)\geq s\right)\mathrm{d}s+\int_{N^{-8/5}}^{\pi}\mathbb{P}\left(\min_{F\in\mathcal{F}(P_N)}\mathrm{vol}_{n-1}(F)\geq s\right)\mathrm{d}s.\end{aligned}$$

By (84) with a new constant $c_3$

$$\mathbb{E}\left(\min_{F\in\mathcal{F}(P_N)}\mathrm{vol}_{n-1}(F)\right)\leq N^{-8/5}+\frac{c_3}{N^{8/3}}\int_{N^{-8/5}}^{\pi}s^{-5/3}\,\mathrm{d}s\leq N^{-8/5}+\frac{3}{2}c_3N^{-8/5}.\qquad\square$$

# 8 Expected minimum lower bound, $n\geq 3$

Now we prove the estimates from below in Theorem 2, $n\geq 3$. In order to do this we estimate the probability that there is a facet of the random polytope $[\xi_1,\dots,\xi_N]$ that has an $(n-1)$-dimensional volume less than $t$ from above. Proposition 16 assures that this probability is less than $c_2t^2N^3$ when $n\geq 4$, and it is less than $c_2t^{5/3}N^{8/5}$ when $n=3$. Consequently, the probability that this volume is greater than $t$ is greater than $1-c_2t^2N^3$ in the case of dimension greater or equal to 4, and greater than $1-c_2t^{\frac{5}{3}}N^{\frac{8}{3}}$ in dimension 3. Then we choose $t=N^{-\frac{3}{2}}$ and $t=N^{-\frac{8}{5}}$.

**Proposition 16.** *Let $\xi_1,\dots,\xi_N$ be i.i.d. random points chosen from $S^{n-1}$.*
*(i) For $n=3$,*

$$\mathbb{P}\big\{\exists i_1,i_2,i_3:[\xi_{i_1},\xi_{i_2},\xi_{i_3}]\in\mathcal{F}([\xi_1,\dots,\xi_N])\text{ and }\mathrm{vol}_2([\xi_{i_1},\xi_{i_2},\xi_{i_3}])\leq t\big\}\leq 57\pi N^{\frac{8}{3}}t^{\frac{5}{3}} \tag{85}$$

*(ii) For $n\geq 4$ and $b_n$ given by (39) and for all sufficiently big $N$,*

$$\begin{aligned}&\mathbb{P}\big\{\exists i_1,\dots,i_n:[\xi_{i_1},\dots,\xi_{i_n}]\in\mathcal{F}([\xi_1,\dots,\xi_N])\text{ and }\mathrm{vol}_{n-1}([\xi_{i_1},\dots,\xi_{i_n}])\leq t\big\}\\&\leq 4\sqrt{2}b_{n-1}t^2N^3\frac{|S^{n-2}|^3}{n|S^{n-1}|^2}(n-1)^{n-4}(n-4)!.\end{aligned} \tag{86}$$

In fact, the estimates in Proposition 16 give the optimal orders in $t$ and $N$. We skip the proofs of the lower estimates because we do not use them for the proofs of our theorems. The proof of Proposition 16 appears after the proof of the following lemma.

**Lemma 17.** *For $n \in \mathbb{N}$ with $n \geq 4$ there is $N_n \in \mathbb{N}$ such that for all $N \in \mathbb{N}$ with $N \geq N_n$,*

$$\begin{aligned} &\int_0^1 (1-p^2)^{\frac{n^2-4n+1}{2}} \left(1 - \frac{|S^{n-2}|}{|S^{n-1}|} \int_p^1 (1-s^2)^{\frac{n-3}{2}} \,\mathrm{d}s\right)^{N-n} \mathrm{d}p \\ &\qquad \leq 2\sqrt{2}(n-1)^{n-4} \left(\frac{|S^{n-1}|}{|S^{n-2}|}\right)^{n-3} \frac{(n-4)!}{(N-n)^{n-3}}. \end{aligned} \tag{87}$$

*Proof.* Before we begin with the actual proof we want to make sure that the left hand side expression in (87) is well defined, i.e.

$$\frac{|S^{n-2}|}{|S^{n-1}|} \int_0^1 (1-s^2)^{\frac{n-3}{2}} \,\mathrm{d}s \leq 1.$$

Indeed,

$$\frac{|S^{n-2}|}{|S^{n-1}|} \int_0^1 (1-s^2)^{\frac{n-3}{2}} \,\mathrm{d}s = \frac{1}{2} \frac{n-1}{n\sqrt{\pi}} \frac{\Gamma(\frac{n+2}{2})}{\Gamma(\frac{n+1}{2})} \frac{\Gamma(\frac{n-1}{2})\Gamma(\frac{1}{2})}{\Gamma(\frac{n}{2})} = \frac{1}{2}.$$

Now we show (87). For $p \in [0,1]$ we have

$$\int_p^1 (1-s^2)^{\frac{n-3}{2}} \,\mathrm{d}s \geq \int_p^1 s(1-s^2)^{\frac{n-3}{2}} \,\mathrm{d}s = \frac{1}{n-1} \left[-(1-s^2)^{\frac{n-1}{2}}\right]_p^1 = \frac{(1-p^2)^{\frac{n-1}{2}}}{n-1}. \tag{88}$$

We note that

$$\frac{1}{n-1} \frac{|S^{n-2}|}{|S^{n-1}|} \leq \frac{1}{2\sqrt{\pi}}. \tag{89}$$

Indeed, since $\Gamma(z+1) = z\Gamma(z)$ and since $\Gamma(z)$ is increasing for $z \geq 2$

$$\frac{1}{n-1} \frac{|S^{n-2}|}{|S^{n-1}|} = \frac{1}{n\sqrt{\pi}} \frac{\Gamma(\frac{n}{2}+1)}{\Gamma(\frac{n-1}{2}+1)} = \frac{1}{2\sqrt{\pi}} \frac{\Gamma(\frac{n}{2})}{\Gamma(\frac{n+1}{2})} \leq \frac{1}{2\sqrt{\pi}}.$$

By (88), and with the substitution $s = (1-p^2)^{\frac{n-1}{2}}$, i.e. $p = \left(1 - s^{\frac{2}{n-1}}\right)^{1/2}$ and $\mathrm{d}p = -\frac{s^{-\frac{n-3}{n-1}}}{(n-1)\sqrt{1-s^{\frac{2}{n-1}}}} \,\mathrm{d}s$, the left hand side of (87) is less than

$$\begin{aligned} &\int_0^1 (1-p^2)^{\frac{n^2-4n+1}{2}} \left(1 - \frac{|S^{n-2}|}{|S^{n-1}|} \frac{(1-p^2)^{\frac{n-1}{2}}}{n-1}\right)^{N-n} \mathrm{d}p \\ &= \frac{1}{n-1} \int_0^1 s^{n-4} \left(1 - \frac{|S^{n-2}|}{|S^{n-1}|} \frac{s}{n-1}\right)^{N-n} \left(1 - s^{\frac{2}{n-1}}\right)^{-1/2} \mathrm{d}s \\ &= \frac{1}{n-1} \int_0^{2^{-\frac{n-1}{2}}} s^{n-4} \left(1 - \frac{|S^{n-2}|}{|S^{n-1}|} \frac{s}{n-1}\right)^{N-n} \left(1 - s^{\frac{2}{n-1}}\right)^{-1/2} \mathrm{d}s \\ &\quad + \frac{1}{n-1} \int_{2^{-\frac{n-1}{2}}}^1 s^{n-4} \left(1 - \frac{|S^{n-2}|}{|S^{n-1}|} \frac{s}{n-1}\right)^{N-n} \left(1 - s^{\frac{2}{n-1}}\right)^{-1/2} \mathrm{d}s. \end{aligned} \tag{90}$$

We estimate the first summand. Since $\left(1-s^{\frac{2}{n-1}}\right)^{1/2} \geq \frac{1}{\sqrt{2}}$ for $s \in [0, 2^{-\frac{n-1}{2}}]$, the first summand is smaller than

$$\frac{\sqrt{2}}{n-1}\int_0^{2^{-\frac{n-1}{2}}} s^{n-4}\left(1-\frac{|S^{n-2}|}{|S^{n-1}|}\frac{s}{n-1}\right)^{N-n} \mathrm{d}s.$$

We substitute $t=\frac{|S^{n-2}|}{|S^{n-1}|}\frac{s}{n-1}$ and $\mathrm{d}s=(n-1)\frac{|S^{n-1}|}{|S^{n-2}|}\,\mathrm{d}t$ and get

$$\frac{\sqrt{2}}{n-1}\left(\frac{|S^{n-1}|}{|S^{n-2}|}(n-1)\right)^{n-3}\int_0^{\frac{2^{-\frac{n-1}{2}}|S^{n-2}|}{(n-1)|S^{n-1}|}} t^{n-4}\left(1-t\right)^{N-n}\mathrm{d}t. \tag{91}$$

By (89) we get $\frac{2^{-\frac{n-1}{2}}|S^{n-2}|}{(n-1)|S^{n-1}|} \leq 1$. Therefore the integral in (91) is smaller than

$$\int_0^1 t^{n-4}\left(1-t\right)^{N-n}\mathrm{d}t = B(n-3, N-n+1) \leq \frac{\Gamma(n-3)}{(N-n)^{n-3}}.$$

Now we estimate the second summand of (90). In the integral we replace $s^{n-4}$ by $s^{-\frac{n-3}{n-1}}$. Since $s^{n-4} \leq s^{-\frac{n-3}{n-1}}$ the integral does not decrease. Therefore, the second summand of (90) is smaller than

$$\begin{aligned}
&\frac{1}{n-1}\int_{2^{-\frac{n-1}{2}}}^1 \exp\left(-\frac{(N-n)|S^{n-2}|}{(n-1)|S^{n-1}|}s\right) s^{-\frac{n-3}{n-1}}\left(1-s^{\frac{2}{n-1}}\right)^{-1/2}\mathrm{d}s \\
&\quad\leq \frac{1}{n-1}\exp\left(-\frac{(N-n)|S^{n-2}|}{(n-1)|S^{n-1}|}2^{-\frac{n-1}{2}}\right)\int_{2^{-\frac{n-1}{2}}}^1 s^{-\frac{n-3}{n-1}}\left(1-s^{\frac{2}{n-1}}\right)^{-1/2}\mathrm{d}s \\
&\quad= \frac{1}{\sqrt{2}}\exp\left(-\frac{(N-n)|S^{n-2}|}{2^{\frac{n-1}{2}}(n-1)|S^{n-1}|}\right).
\end{aligned}$$

Altogether (87) is less than

$$\frac{\sqrt{2}}{n-1}\left(\frac{(n-1)|S^{n-1}|}{|S^{n-2}|}\right)^{n-3}\frac{(n-4)!}{(N-n)^{n-3}}+\frac{1}{\sqrt{2}}\exp\left(-\frac{(N-n)|S^{n-2}|}{2^{-\frac{n-1}{2}}(n-1)|S^{n-1}|}\right). \tag{92}$$

For each $n \geq 4$ there is $N_n$ such that for all $N \geq N_n$ the first summand in (92) is larger than the second. It follows that for each $n \geq 4$ there is $N_n$ such that for all $N \geq N_n$ we have that (87) is less than

$$\frac{2\sqrt{2}}{n-1}\left(\frac{(n-1)|S^{n-1}|}{|S^{n-2}|}\right)^{n-3}\frac{(n-4)!}{(N-n)^{n-3}}. \qquad \square$$

**Proof of Proposition 16**. *(ii)* We have

$$
\begin{aligned}
&\mathbb{P}\big\{\exists i_1,\dots,i_n : [\xi_{i_1},\dots,\xi_{i_n}] \in \mathcal{F}([\xi_1,\dots,\xi_N]) \text{ and } \mathrm{vol}_{n-1}([\xi_{i_1},\dots,\xi_{i_n}]) \le t\big\} \\
&\quad \le \sum_{i_1,\dots,i_n} \mathbb{P}\big\{[\xi_{i_1},\dots,\xi_{i_n}] \in \mathcal{F}([\xi_1,\dots,\xi_N]) \text{ and } \mathrm{vol}_{n-1}([\xi_{i_1},\dots,\xi_{i_n}]) \le t\big\} \\
&\quad = \frac{\binom{N}{n}}{|S^{n-1}|^n} \int_{(S^{n-1})^n} \chi\big(\mathrm{vol}_{n-1}([\xi_1,\dots,\xi_n]) \le t\big) \left( \left( \frac{\mathrm{vol}_{n-1}(\mathrm{aff}(\xi_{i_1},\dots,\xi_{i_n})^+ \cap S^{n-1})}{|S^{n-1}|} \right)^{N-n} \right. \\
&\qquad\qquad \left. + \left( \frac{\mathrm{vol}_{n-1}(\mathrm{aff}(\xi_{i_1},\dots,\xi_{i_n})^- \cap S^{n-1})}{|S^{n-1}|} \right)^{N-n} \right) \mathrm{d}\xi_1 \cdots \mathrm{d}\xi_n, \qquad (93)
\end{aligned}
$$

where $\mathrm{aff}(\xi_{i_1},\dots,\xi_{i_n})^+$ is the halfspace containing the origin. The case where both halfspaces contain the origin does not play any role in the computation. Clearly,

$$\mathrm{vol}_{n-1}\big(\mathrm{aff}(\xi_{i_1},\dots,\xi_{i_n})^- \cap S^{n-1}\big) \le |S^{n-1}|/2$$

and consequently,

$$
\begin{aligned}
&\frac{\binom{N}{n}}{|S^{n-1}|^n} \int_{(S^{n-1})^n} \chi\big(\mathrm{vol}_{n-1}([\xi_1,\dots,\xi_n]) \le t\big) \left( \frac{\mathrm{vol}_{n-1}\big(\mathrm{aff}(\xi_{i_1},\dots,\xi_{i_n})^- \cap S^{n-1}\big)}{|S^{n-1}|} \right)^{N-n} \mathrm{d}\xi_1 \cdots \mathrm{d}\xi_n \\
&\qquad \le 2^{-N+n} \binom{N}{n} \le 2^{-N} N^n. \qquad (94)
\end{aligned}
$$

By the Blaschke-Petkantschin formula (Lemma 7), the other summand of (93) equals

$$
\begin{aligned}
&\binom{N}{n} \frac{(n-1)!}{|S^{n-1}|^n} \int_{S^{n-1}} \int_0^1 \int_{(H(\theta,p)\cap S^{n-1})^n} \left( \frac{\mathrm{vol}_{n-1}\big(H(\theta,p)^+ \cap S^{n-1}\big)}{|S^{n-1}|} \right)^{N-n} \frac{1}{(1-p^2)^{\frac{n}{2}}} \qquad (95) \\
&\qquad \chi\big(\mathrm{vol}_{n-1}([\eta_1,\dots,\eta_n]) \le t\big)\, \mathrm{vol}_{n-1}([\eta_1,\dots,\eta_n])\, \mathrm{d}\eta_1 \cdots \mathrm{d}\eta_n\, \mathrm{d}p\, \mathrm{d}\theta.
\end{aligned}
$$

By (1) in Lemma 3 for $n \ge 3$,

$$\frac{\mathrm{vol}_{n-1}(H(\theta,p)^+ \cap S^{n-1})}{|S^{n-1}|} = 1 - \frac{|S^{n-2}|}{|S^{n-1}|} \int_p^1 (1-s^2)^{\frac{n-3}{2}}\, \mathrm{d}s.$$

The integral in (95) is rotationally invariant with respect to $\theta$. By this and (94),

$$
\begin{aligned}
&\mathbb{P}\big\{\exists i_1,\dots,i_n : [\xi_{i_1},\dots,\xi_{i_n}] \in \mathcal{F}([\xi_1,\dots,\xi_N]) \text{ and } \mathrm{vol}_{n-1}([\xi_{i_1},\dots,\xi_{i_n}]) \le t\big\} \\
&\quad \le \frac{N^n}{2^N} + \binom{N}{n} \frac{(n-1)!}{|S^{n-1}|^{n-1}} \int_0^1 \int_{(H(e_1,p)\cap S^{n-1})^n} \left( 1 - \frac{|S^{n-2}|}{|S^{n-1}|} \int_p^1 (1-s^2)^{\frac{n-3}{2}}\, \mathrm{d}s \right)^{N-n} \frac{1}{(1-p^2)^{\frac{n}{2}}} \\
&\qquad\qquad \chi\big(\mathrm{vol}_{n-1}([\eta_1,\dots,\eta_n]) \le t\big)\, \mathrm{vol}_{n-1}\big([\eta_1,\dots,\eta_n]\big)\, \mathrm{d}\eta_1 \cdots \mathrm{d}\eta_n\, \mathrm{d}p. \qquad (96)
\end{aligned}
$$

With $\eta_i = \sqrt{1-p^2}\zeta_i$, $\mathrm{d}\eta_i = (1-p^2)^{\frac{n-2}{2}}\, \mathrm{d}\zeta_i$,

$$
\begin{aligned}
&\int_{(H(e_1,p)\cap S^{n-1})^n} \chi\big(\mathrm{vol}_{n-1}([\eta_1,\dots,\eta_n]) \le t\big)\, \mathrm{vol}_{n-1}([\eta_1,\dots,\eta_n])\, \mathrm{d}\eta_1 \cdots \mathrm{d}\eta_n \\
&= \int_{(S^{n-2})^n} \chi\left( \mathrm{vol}_{n-1}([\zeta_1,\dots,\zeta_n]) \le \frac{t}{(1-p^2)^{\frac{n-1}{2}}} \right) \mathrm{vol}_{n-1}([\zeta_1,\dots,\zeta_n]) (1-p^2)^{\frac{n^2-n-1}{2}}\, \mathrm{d}\zeta_1 \cdots \mathrm{d}\zeta_n.
\end{aligned}
$$

By (64),

$$\int_{(S^{n-2})^n} \chi\left(\mathrm{vol}_{n-1}([\zeta_1,\dots,\zeta_n]) \le t(1-p^2)^{-\frac{n-1}{2}}\right) \mathrm{vol}_{n-1}([\zeta_1,\dots,\zeta_n])\,\mathrm{d}\zeta_1\cdots\mathrm{d}\zeta_n$$

$$= |S^{n-2}|^n \int_0^{t(1-p^2)^{-\frac{n-1}{2}}} \mathbb{P}\left(s \le \mathrm{vol}_{n-1}([\zeta_1,\dots,\zeta_n]) \le t(1-p^2)^{-\frac{n-1}{2}}\right)\mathrm{d}s$$

$$\le |S^{n-2}|^n t(1-p^2)^{-\frac{n-1}{2}} \cdot \mathbb{P}\left(\mathrm{vol}_{n-1}([\zeta_1,\dots,\zeta_n]) \le t(1-p^2)^{-\frac{n-1}{2}}\right).$$

By (37) the previous expression is smaller than

$$\int_{(H(e_1,p)\cap S^{n-1})^n} \chi\big(\mathrm{vol}_{n-1}([\eta_1,\dots,\eta_n]) \le t\big)\,\mathrm{vol}_{n-1}([\eta_1,\dots,\eta_n])\,\mathrm{d}\eta_1\cdots\mathrm{d}\eta_n$$

$$\le b_{n-1}|S^{n-2}|^n \left(t(1-p^2)^{-\frac{n-1}{2}}\right)^2 (1-p^2)^{\frac{n^2-n-1}{2}} = b_{n-1}|S^{n-2}|^n t^2 (1-p^2)^{\frac{n^2-3n+1}{2}},$$

where $b_{n-1}$ is given by (39). This in (96) gives

$$\mathbb{P}\big\{\exists i_1,\dots,i_n : [\xi_{i_1},\dots,\xi_{i_n}] \in \mathcal{F}([\xi_1,\dots,\xi_N]) \text{ and } \mathrm{vol}_{n-1}([\xi_{i_1},\dots,\xi_{i_n}]) \le t\big\}$$
$$\le \frac{N^n}{2^N} + b_{n-1}t^2\binom{N}{n}\frac{(n-1)!|S^{n-2}|^n}{|S^{n-1}|^{n-1}}\int_0^1\left(1-\frac{|S^{n-2}|}{|S^{n-1}|}\int_p^1(1-s^2)^{\frac{n-3}{2}}\,\mathrm{d}s\right)^{N-n}(1-p^2)^{\frac{n^2-4n+1}{2}}\,\mathrm{d}p.$$

Now we apply Lemma 17 (see also [16, Lemma 4.10]). For every $n$ there is $N_n$ such that for all $N$ with $N \ge N_n$ the above expression is smaller than

$$\frac{N^n}{2^N} + b_{n-1}t^2\binom{N}{n}\frac{(n-1)!|S^{n-2}|^n}{|S^{n-1}|^{n-1}}2\sqrt{2}(n-1)^{n-4}\left(\frac{|S^{n-1}|}{|S^{n-2}|}\right)^{n-3}\frac{(n-4)!}{(N-n)^{n-3}}$$
$$\le \frac{N^n}{2^N} + 2\sqrt{2}b_{n-1}t^2\frac{N^n}{(N-n)^{n-3}}\frac{|S^{n-2}|^3}{n|S^{n-1}|^2}(n-1)^{n-4}(n-4)!.$$

For large enough $N$ this expression is smaller than

$$4\sqrt{2}b_{n-1}t^2N^3\frac{|S^{n-2}|^3}{n|S^{n-1}|^2}(n-1)^{n-4}(n-4)!.$$

Thus we have shown (86).

*(i)* We consider now the case $n=3$. By (96),

$$\mathbb{P}\big\{\exists i_1,i_2,i_3 : [\xi_{i_1},\xi_{i_2},\xi_{i_3}] \in \mathcal{F}([\xi_1,\dots,\xi_N]) \text{ and } \mathrm{vol}_2([\xi_{i_1},\xi_{i_2},\xi_{i_3}]) \le t\big\}$$
$$\le \frac{N^3}{2^N} + \binom{N}{3}\frac{2}{|S^2|^2}\int_0^1\int_{(H(e_1,p)\cap S^2)^3}\left(1-\frac{|S^1|}{|S^2|}\int_p^1\mathrm{d}s\right)^{N-3}\frac{\chi\big(\mathrm{vol}_2([\eta_1,\eta_2,\eta_3]) \le t\big)}{(1-p^2)^{3/2}}$$
$$\mathrm{vol}_2([\eta_1,\eta_2,\eta_3])\,\mathrm{d}\eta_1\,\mathrm{d}\eta_2\,\mathrm{d}\eta_3\,\mathrm{d}p$$
$$\le \frac{N^3}{2^N} + \frac{N^3}{48\pi^2}\int_0^1\int_{(H(e_1,p)\cap S^2)^3}\left(\frac{1+p}{2}\right)^{N-3}\frac{\chi\big(\mathrm{vol}_2([\eta_1,\eta_2,\eta_3]) \le t\big)}{(1-p^2)^{3/2}}$$
$$\mathrm{vol}_2([\eta_1,\eta_2,\eta_3])\,\mathrm{d}\eta_1\,\mathrm{d}\eta_2\,\mathrm{d}\eta_3\,\mathrm{d}p.$$

With $\eta_i=\sqrt{1-p^2}\zeta_i$, $\mathrm{d}\eta_i=\sqrt{1-p^2}\,\mathrm{d}\zeta_i$,

$$\begin{aligned}&\int_{(H(e_1,p)\cap S^2)^3}\chi\big(\mathrm{vol}_2([\eta_1,\eta_2,\eta_3])\leq t\big)\,\mathrm{vol}_2([\eta_1,\eta_2,\eta_3])\,\mathrm{d}\eta_1\,\mathrm{d}\eta_2\,\mathrm{d}\eta_3\\&=\int_{(S^1)^3}\chi\left(\mathrm{vol}_2([\zeta_1,\zeta_2,\zeta_3])\leq\frac{t}{1-p^2}\right)\mathrm{vol}_2([\zeta_1,\zeta_2,\zeta_3])(1-p^2)^{5/2}\,\mathrm{d}\zeta_1\,\mathrm{d}\zeta_2\,\mathrm{d}\zeta_3.\end{aligned}$$

Again, by (64) this equals

$$\begin{aligned}&(1-p^2)^{5/2}|S^1|^3\int_0^{\frac{t}{1-p^2}}\mathbb{P}\left(s\leq\mathrm{vol}_2([\zeta_1,\zeta_2,\zeta_3])\leq\frac{t}{1-p^2}\right)\mathrm{d}s\\&\qquad\leq(1-p^2)^{3/2}|S^1|^3t\cdot\mathbb{P}\left(\mathrm{vol}_2([\zeta_1,\zeta_2,\zeta_3])\leq\frac{t}{1-p^2}\right).\end{aligned}$$

By (36) this is smaller than

$$(1-p^2)^{3/2}|S^1|^3t\cdot342\cdot\left(\frac{t}{1-p^2}\right)^{2/3}=342\cdot8\pi^3t^{\frac{5}{3}}(1-p^2)^{\frac{5}{6}}.$$

Therefore

$$\begin{aligned}&\mathbb{P}\left\{\exists i_1,i_2,i_3:[\xi_{i_1},\xi_{i_2},\xi_{i_3}]\in\mathcal{F}([\xi_1,\dots,\xi_N])\text{ and }\mathrm{vol}_2([\xi_{i_1},\xi_{i_2},\xi_{i_3}])\leq t\right\}\\&\leq\frac{N^3}{2^N}+\frac{N^3}{48\pi^2}342\cdot8\pi^3t^{\frac{5}{3}}\int_0^1(1-p^2)^{-\frac{3}{2}}\left(\frac{1+p}{2}\right)^{N-3}(1-p^2)^{\frac{5}{6}}\,\mathrm{d}p\\&=\frac{N^3}{2^N}+N^3t^{\frac{5}{3}}57\pi\int_0^1(1-p^2)^{-\frac{2}{3}}\left(\frac{1+p}{2}\right)^{N-3}\mathrm{d}p.\end{aligned}$$

We put $u=\frac{1+p}{2}$, i.e. $p=2u-1$. The previous expression equals

$$\begin{aligned}&\frac{N^3}{2^N}+N^3t^{\frac{5}{3}}2\cdot57\pi\int_{1/2}^1(4u-4u^2)^{-\frac{2}{3}}u^{N-3}\,\mathrm{d}u=\frac{N^3}{2^N}+N^3t^{\frac{5}{3}}\frac{57\pi}{2^{\frac{1}{3}}}\int_{1/2}^1(1-u)^{-\frac{2}{3}}u^{N-\frac{11}{3}}\,\mathrm{d}u\\&\leq\frac{N^3}{2^N}+N^3t^{\frac{5}{3}}\frac{57\pi}{2^{\frac{1}{3}}}B\left(\frac{1}{3},N-\frac{8}{3}\right)=\frac{N^3}{2^N}+N^3t^{\frac{5}{3}}\frac{57\pi}{2^{\frac{1}{3}}}\frac{\Gamma(\frac{1}{3})\Gamma(N-\frac{8}{3})}{\Gamma(N-\frac{7}{3})}.\end{aligned}$$

Since $\lim_{x\to\infty}\frac{\Gamma(x+\alpha)}{x^\alpha\Gamma(x)}=1$, for large enough $N$ we get, with a bigger constant,

$$\mathbb{P}\Big\{\exists i_1,i_2,i_3:[\xi_{i_1},\xi_{i_2},\xi_{i_3}]\in\mathcal{F}([\xi_1,\dots,\xi_N])\text{ and }\mathrm{vol}_2([\xi_{i_1},\xi_{i_2},\xi_{i_3}])\leq t\Big\}\leq57\pi\cdot N^{\frac{8}{3}}t^{\frac{5}{3}}.\square$$

**Proof of the lower estimates in Theorem 2,** $n\geq3$. (iii) This is the case $n\geq4$. Let $P_N=[\xi_1,\dots,\xi_N]$ be the random polytope. We have $\min_{F\in\mathcal{F}(P_N)}\mathrm{vol}_{n-1}(F)\geq t$ if and only if

$$\forall i_1,\dots,i_n\in\{1,\dots,N\}:[\xi_{i_1},\dots,\xi_{i_n}]\notin\mathcal{F}(P_N)\text{ or }\mathrm{vol}_{n-1}([\xi_{i_1},\dots,\xi_{i_n}])\geq t.$$

The measure of this event is

$$\begin{aligned}&\mathbb{P}\Big\{\forall i_1,\dots,i_n\in\{1,\dots,N\}:[\xi_{i_1},\dots,\xi_{i_n}]\notin\mathcal{F}(P_N)\text{ or }\mathrm{vol}_{n-1}([\xi_{i_1},\dots,\xi_{i_n}])\geq t\Big\}\\&\quad=1-\mathbb{P}\Big\{\exists i_1,\dots,i_n:[\xi_{i_1},\dots,\xi_{i_n}]\in\mathcal{F}(P_N)\text{ and }\mathrm{vol}_{n-1}([\xi_{i_1},\dots,\xi_{i_n}])<t\Big\}.\end{aligned}$$

By (86) the previous expression is greater than

$$1-4\sqrt{2}b_{n-1}t^2N^3\frac{|S^{n-2}|^3}{n|S^{n-1}|^2}(n-1)^{n-4}(n-4)!,$$

where $b_{n-1}$ is given by (39). We choose

$$t=\left(8\sqrt{2}b_{n-1}N^3\frac{|S^{n-2}|^3}{n|S^{n-1}|^2}(n-1)^{n-4}(n-4)!\right)^{-\frac{1}{2}}$$

and we get

$$\mathbb{P}\big\{\forall i_1,\dots,i_n:[\xi_{i_1},\dots,\xi_{i_n}]\notin\mathcal{F}(P_N)\text{ or }\operatorname{vol}_{n-1}([\xi_{i_1},\dots,\xi_{i_n}])\geq t\big\}\geq 1/2.$$

By Markov's inequality,

$$\mathbb{E}\min_{F\in\mathcal{F}(P_N)}\operatorname{vol}_{n-1}(F)\geq\frac{t}{2}=\frac{N^{-3/2}}{2}\left(8\sqrt{2}b_{n-1}\frac{|S^{n-2}|^3}{n|S^{n-1}|^2}(n-1)^{n-4}(n-4)!\right)^{-1/2}.$$

(ii) In the case of dimension 3, we proceed similarly to the previous case (case (iii), $n\geq 4$). We apply (85), set $t=(\frac{1}{114\pi})^{\frac{3}{5}}N^{-8/5}$ and apply Markov's inequality to get

$$\mathbb{E}\min_{F\in\mathcal{F}(P_N)}\operatorname{vol}_2(F)\geq\frac{t}{2}=\frac{1}{2}\left(\frac{1}{114\pi}\right)^{3/5}N^{-8/5}.\qquad\square$$

# 9 Expected minimum, $n=2$

**Proof of Theorem 2, $n=2$.** (i) This is known and it can be deduced from the expected length of the smallest gap in a sample of uniformly random points in $S^1$, which in turn follows from the expected length of the smallest gap in a sample of uniformly random points in $[0,1]$. From [11, pages 63–64] (based on an argument from [13]) we have that for a sample of $N$ uniformly random points from $S^1$ the expectation of the minimum geodesic distance among pairs of points is $2\pi/N^2$. To relate the arc length with the edge length of the convex hull of the sample, it is enough to notice that for $N\geq 3$ the length of the smallest arc is at most $2\pi/3$ and therefore the edge length is at most the arc length and at least $\frac{3\sqrt{3}}{2\pi}$ times the arc length. The claim follows. □